\documentclass{amsproc}
\usepackage{graphicx}

\newtheorem{theorem}{Theorem}[section]
\newtheorem{lemma}[theorem]{Lemma}
\newtheorem{proposition}[theorem]{Proposition}

\newtheorem{definition}[theorem]{Definition}

\begin{document}

\title[Navier-Stokes equations generate all  hyperbolic dynamics]{Navier-Stokes equations under Marangoni boundary conditions generate all  hyperbolic dynamics}

\author{Sergei Vakulenko}

\address{Institute for Mechanical Engineering Problems,
 Bolshoy pr. V. O. 61, 199178, Saint Petersburg,   and
Saint Petersburg National Research University of Information Technologies, Mechanics and Optics, 197101, St. Petersburg, Russia}


\begin{abstract}
In this paper, we consider
  dynamics defined by the Navier-Stokes equations under the Marangoni boundary condtions in a two dimensional domain.  This model of fluid dynamics  involve fundamental physical effects: convection, diffusion and capillary forces.
The main result  is as follows:
 local semiflows,  defined by the corresponding initial boundary value problem, can  generate   all possible
 structurally stable  dynamics defined by $C^1$ smooth vector fields on  compact smooth
 manifolds  (up to an orbital topological equivalence).  To generate a prescribed dynamics, it is sufficient to 
  adjust some parameters in the  equations, namely, the Prandtl number and an external heat source.

\end{abstract}

\maketitle

\section{Introduction}
\label{intro}


In this paper, we state an analytical proof of 
 existence of strange attractors for a model of  fluid dynamics.  The hypothesis that  turbulence can be  connected with strange (chaotic) attractors was pionereed in papers \cite{RT, NRT}.  
We consider the initial boundary values problem (IBVP) defined by
the Navier-Stokes (NS) equations   in the two-dimensional case
  with the Marangoni boundary conditions. These equations
 present a   model
for non-compressible fluid dynamics, which involves fundamental physical effects: convection,  heat transfer
and capillarity. This model describes surface driven the B\'enard-Marangoni  convection  leading to interesting phenomena, for example,
B\'enard cells, 
(see \cite{Bracard,  Drazin, GJ}  and references therein).

The main result of this paper can be outlined as follows:
 local semiflows, induced by the NS equations with the Marangoni boundary conditions, can  generate   all possible
 hyperbolic  dynamics defined by $C^1$-smooth vector fields on  finite dimensional compact smooth
 manifolds  (up to an orbital topological equivalence).  To generate a prescribed dynamics, it is sufficient to
 adjust some parameters in the NS equations, namely, the Prandtl number and a heat source in the heat transfer
 equation.  The well known examples of  hyperbolic dynamics
with a ''chaotic" behaviour are  the Anosov flows, the Smale A-axiom systems and  the Smale horseshoes  \cite{Ru, Katok, Holmes}.

Although numerous works were dedicated to the NS equations,
( for example, \cite{Temam, Tem, Lions1,Lions2}),
nonetheless, up to now
   general results  on  dynamical complexity of dissipative dynamical systems defined by
the NS equations do not exist.   Results on chaos existence were obtained 
 for quasilinear parabolic equations and reaction-diffusion systems  \cite{Pol1, Pol2, Vak4} by a special approach.
The main ingredient of this  approach is the method of realization of
vector fields (RVF) proposed by P. Pol\'a\v cik \cite{Pol1, Pol2}.
  Let us outline  briefly the RVF  method and some
 previously obtained results.

Let us consider an  initial boundary value problem associated with a system of PDE's. Assume this problem
 involves some parameters $\mathcal P$  (for example, in the Marangoni problem  the parameters are a  heat source and the Prandtl
number).  We obtain a family of local semiflows $S^t_{\mathcal P}$ generated by these initial boundary value problems,
where each semiflow depends on the corresponding parameter value  ${\mathcal P}$.
Suppose
for an integer $n > 0$, there is an appropriate choice of the parameter ${\mathcal P}_n$
such that the corresponding initial boundary value problem
generates a global semiflow $S^t({\mathcal P}_n)$  possessing an $n$-
dimensional finite $C^{1}$ -smooth positively invariant manifold  ${\mathcal M}_n$ (we can suppose, for simplicity, that this manifold is
diffeomorphic to the unit  ball ${\mathcal B}^n \subset {\bf R}^n$).
The semiflow $S^t({\mathcal P}_n)$, restricted to ${\mathcal M}_n$   (a local inertial form), is
 defined  by a $C^{1}$-vector field $Q$ on ${\mathcal B}^n$. Then we say that the family $S^t_{\mathcal P}$
realizes the vector field $Q$. We say that this family $\epsilon$- realizes a vector field $\bar Q$
if the field $Q$ is a $\epsilon$ -perturbation of the field $\bar Q$ in $C^1({\mathcal B}^n)$ (we consider $C^1$-norms
in order to apply the theorem on persistence of hyperbolic  sets \cite{Ru,Katok}).

By the RVF method,
it has been shown that semiflows
associated with some special quasilinear parabolic equations in two dimensional
domains can generate complicated hyperbolic  sets \cite{Pol2,DaP}.
For reaction-diffusion systems of a special form the RVF method allows to prove 
existence
of chaotic regimes \cite{Vak4}.
One can show that, for any integer $n$, semiflows induced by these systems  can realize
a dense set of the fields in the space of all $C^1$-smooth vector fields on ${\mathcal B}^n$ \cite{Vak4}. Therefore, such systems
generate any structurally stable
 ( persistent under
sufficiently small $C^1$- perturbations) dynamics, up to orbital topological equivalence  \cite{Ru, Katok}.
The corresponding families of the semiflows  can be called {\em maximally dynamically complex}.
If a family of semiflows enjoys this dynamical complexity property, this family generates all
 compact invariant hyperbolic dynamics  on finite dimensional compact smooth manifolds.

By this terminology, the main result of this paper is as follows.
The family of semiflows associated with our IBVP possesses the property
of maximal dynamical complexity (see Theorem \ref{maint}). 

Using the RVF  we encounter  two main technical difficulties.
First,  in previous papers \cite{Pol1, Pol2, DaP, Vak4}
the RVF has been used only for systems with non-polynomial nonlinearities whereas the NS equations involve quadratic ones.
In the first part of the paper we
 overcome this difficulty. We consider systems of differential equations with quadratic nonlinearities
\begin{equation} \label{0.1}
 \frac{dX}{dt}={ K}(X,X) + { M}X + { f},\quad X \in {\bf R}^N, 
\end{equation}
where  $  X(t) $ is a unknown function, $X=(X_1, ..., X_N) \in {\bf R}^N$,  $K(X, X)$ is a bilinear quadratic form, $f \in {\bf R}^N$, $M$ is a $N \times N$ matrix.
 Systems (\ref{0.1}) have important applications, in particular, in chemistry, where
 they describe bimolecular chemical reactions \cite{Zab}, and for population dynamics. 
Results  on existence of complex dynamics for (\ref{0.1}) were first obtained in \cite{Korz2} (see
 also \cite{Zab}). In the paper \cite{Kynion} it is shown
that, when we vary $N, { K}, { M}$ and $ f$,  systems (\ref{0.1}) realize all polynomial equations $D^n z=p(z,Dz, ..., D^{n-1}z)$
of $n$-th order for all $n$.  The result \cite{Kynion} is based on purely algebraic methods whereas the work \cite{Korz2} uses
slow manifolds.
 In \cite{Stud} the RVF method  is applied to investigate dynamics of  systems (\ref{0.1}).  It is shown
that systems (\ref{0.1}) generate a maximally complex family of semiflows, where parameters
 are $ N,  M$, the bilinear form $K$ and $f$.  
 
In this paper, we are dealing with a more complicated situation, when quadratic forms $K(X, X)$ can not be considered as free parameters. 
This difficulty is not too hard and it can be overcome  by the methods  of 
 the invariant manifold theory \cite{He,CLu, CFNT, Bates,Tem, Wig, Van} 
that allows us to reduce systems (\ref{0.1}) of large dimension to  analogous systems of smaller dimension, where the corresponding $K(X,X)$
can be considered as free parameters  (see section   \ref{sec:5} for more detail).

The second part of the proof resolves a much more sophisticated problem: how to reduce the dynamics defined by our  IBVP to  systems (\ref{0.1}). 
 Physically, this reduction
 describes dynamics of some "main modes" $X$. We  reduce  the Marangoni problem to system (\ref{0.1}) by locally attracting invariant manifolds.
In order to proceed it, we  choose  the parameters in the NS equations in  a special way.
 It allows us to extract a finite set of "main modes" .
Note that  it is impossible to prove existence of locally attracting invariant or inertial manifolds
for the general NS equations, even for two dimensional case \cite{CFNT}, therefore, we really need this special parameter choice.  We proceed this extraction  using spectral properties of
a linear operator $L$ that determines the linearization of the NS equations. This operator has $N$ zero eigenvalues and
all the rest spectrum of $L$ lies in the negative half plane and it is separated by a positive barrier from the imaginary axis.
 To find  the operator $L$ with such spectral barrier (gap)  is not  easy, and the construction of $L$  is a main
 difficulty in the paper.  

We use  a special 
 choice of a function $\eta(x,y)=U_N(y) + \gamma u_1(x,y)$ that determines the heat source, where 
 $\gamma>0$ is a small parameter.  This parameter defines the nonlinearity,  i.e., the IBVP is weakly nonlinear
for  small $\gamma$.  In $\eta(x,y)$ the argument $x$ is a horizontal space variable,
 $y$ is a vertical axis and the Marangoni conditions hold at the boundary $y=0$. 
The terms $U_N(y)$  and $\gamma u_1$  play different roles.  The operator $L$ is defined by   $U_N$ and it does not depend on $\gamma$.
The function $u_1$ defines the matrix ${ M}$ in equations (\ref{0.1}).  

To simplify analysis, we set periodic boundary conditions along $x$. Then,  since $U_N$ depends only on $y$, we can  separate variables in the spectral problem for the linear operator  $L$ (this method is well known, 
  see  \cite{ RL, Drazin, GJ,Bracard,Sengul}).  Eigenfunctions of $L=L_N$ have the form $\exp(ikx) \theta_k(y)$  with eigenvalues $\lambda_k$, where $k=1,2, ... \ $.   For each $k$ we
  reduce the spectral problem to a nonlinear equation for the spectral parameter $\lambda_k$.
  For a special choice of the function
 $U_N$ this equation can be investigated and we  can check explicitly that $L_N$ has a  spectral barrier (gap), mentioned above.  The investigation of this nontrivial nonlinear equation is the most complicated
part of the paper. 

Note that  the choice of $U_N(y)$ admit a transparent physical interpretation. The function $U_N(y)$  consists of two terms, $U_N=H(y, \mu) + \mu W_N(y,\mu)$, 
where $\mu$ is a second small parameter (independent of $\nu$ and $\gamma$). 
   For small $\mu$ the first term $H_N$ is close
to a step function, where the step is located at the point  $y_0(\mu)$, which  lies at the boundary $y=0$, where the Marangoni condition holds. This means that we have   ''a heat shock"
at the boundary.  The second term $\mu W_N$ is a small  polynomial perturbation of $H$. For $U_N=H$ we have the operator with the eigenvalues $\lambda_k$, which are close to zero for bounded
$k$.  A small polynomial term $\mu W_N$ perturbs the spectrum as follows.  For an especial choice of $\mu W_N$
 we have $\lambda_k=0$,  $k=1,2,...,N$ whereas all others $\lambda_k$    
satisfy $ Re \  \lambda_k  < - \delta(\mu)$, where $\delta>0$.
This situation is  similar to the classical  Rayleigh- B\'enard
convection \cite{RL,  Drazin, GJ}, where $N=1$ and $\lambda_k=0$ for a single wave vector $k=k_0$ only. 
A key difference with respect to this well known case is that our special construction for the Marangoni-B\'enard case produces a number of $k$ with $\lambda_k=0$, and an interaction of the corresponding slow modes can generate a complicated dynamics. 
 
The gap property of the operator $L_N$ allows us to proceed
the reduction of the NS equations to (\ref{0.1})
by a quite routine procedure, which uses the well known results of invariant manifold theory \cite{He, Bates,  CFNT, CLu}.
This procedure shows that $K$ and $M$ depends  on the source $\eta$ in a special way, namely,  the coefficients involved in $K$
depend on  the eigenfunctions of $L=L_N$  whereas $M$ is a linear functional of correction  $u_1$. We show that the  the range of this
functional is  dense in the linear space of all $N \times N$ -matrices.
This fact  allows us to apply   our results  on quadratic  systems (\ref{0.1}) from section \ref{sec:5}
and completes  the proof.

The paper is organized as follows.
In the next section   we formulate  the Marangoni problem.   In Section \ref{sec:3} we state the main result.
Section \ref{sec:4} describes  the RVF method for weakly nonlinear systems.
 In Section \ref{sec:7}   it is shown that the Marangoni problem  is well posed
and defines a local semiflow. In Section \ref{Spectrum}  we investigate
a  linear operator describing a linearization of our IBVP,   and show that this operator has a 
spectral gap.  
In Section \ref{sec: 9}   we   prove existence of the finite dimensional 
 invariant manifold. In Section \ref{sec:5}   we consider quadratic systems (\ref{0.1}).  In Section \ref{sec: 10} we check conditions, which is critically important for
the RVF method. Here we show that, for each fixed $N$,  by a  choice $u_1(x,y)$, we can obtain any prescribed  matrix ${\bf M}$.   The complete algorithm of the  RVF method  uses this fact and it
is stated in Section \ref{sec: 11}.

Below we use the following standard convention: all positive constants,  independent
of small parameters $\epsilon, \gamma ... $, are denoted by $c_i, C_j$. To diminish a formidable number of indices $i,j$,
we shall use sometimes the same indices  assuming that the constants can vary from a line to a line.

\section{Marangoni problem for Navier Stokes equations}
\label{sec:2}
We consider
the Navier Stokes system for an ideal incompressible fluid
\begin{equation}
  {\bf v}_t +  ({\bf v} \cdot \nabla) {\bf v} =\nu \Delta {\bf v} - \nabla p,
\label{OB1}
\end{equation}
\begin{equation}
     \nabla \cdot {\bf v} =0,
\label{div}
\end{equation}
\begin{equation}
   u_t +  ({\bf v} \cdot \nabla) u = \Delta u + \eta,
\label{OB2}
\end{equation}
  where  ${\bf v}=(v_1(x, y, t), v_2(x, y, t))^{tr}$,
$u=u(x,y,t), p=p(x,y, t)$ are   unknown functions
defined on $\Omega \times \{ t \ge 0 \}$,
 $\Omega$ is the strip
$(-\infty, \infty) \times [0,h] \subset {\bf R^2}$.
Here $\bf v$ is the fluid velocity, where $v_1$ and $v_2$ are the normal and tangent velocity components,
                           $\nu$ is
the viscosity coefficient, $p$ is the pressure,
$u$ is the temperature,
$\eta(x,y)$ is a  function describing a distributed
heat source,
${\bf v} \cdot \nabla$ denotes the advection operator $v_1 \frac{\partial}{\partial x} +
v_2 \frac{\partial}{\partial y}$. Notice that since in (\ref{OB2}) the thermal diffusion rate equals $1$, the viscosity coefficient
$\nu$ can be identified with the Prandtl number.

The initial conditions are
\begin{equation}
  {\bf v}(x, y, 0)={\bf v}^0(x, y),
\quad { p}(x, y, 0)={ p}^0(x, y),
\quad u(x,y,0)=u^0(x,y).
\label{inidata}
\end{equation}
Let us
suppose that  the unknown functions
are $2\pi$-periodic  in
$x$:
\begin{equation}
  {\bf v}(x, y, t)={\bf v}(x+ 2\pi, y, t),
\quad { p}(x, y, t)={ p}(x+ 2\pi, y, t),
\label{period1}
\end{equation}
\begin{equation}
 u(x, y, t)= u(x+ 2\pi, y, t),
\label{period2}
\end{equation}
and that  $u^0, p^0, {\bf v}^0$ also are  $2\pi$ -periodic in $x$.
The function $u$ satisfies the Neumann boundary conditions:
\begin{equation}
  u_y(x, y, t)\vert_{y=h} =0, \quad u_y(x, y, t)\vert_{y=0} =0.
\label{boundNeum}
\end{equation}
We assume that the surface $y=h$ is free:
\begin{equation}
  v_2(x, h, t) =0,  \quad  \frac{\partial v_1(x, y, t)}{\partial y}\vert_{y=h}=0.
\label{bound5}
\end{equation}
The Marangoni boundary condition at $y=0$ is defined by a relation
connecting the tangent velocity component and the tangent gradient of the temperature:
\begin{equation}
 { v_1}_y (x,y,t)\vert_{y=0}   = - \gamma_0 u_x(x,0,t),
\label{Maran1}
\end{equation}
where $\gamma_0 >0$ is a coefficient (the Marangoni parameter). We set below $\gamma_0=1$ to simplify formulas. For $v_2$ at $y=0$ one has
\begin{equation}
 { v_2} (x,0,t)   = 0.
\label{Maran1h}
\end{equation}
Let us assume that
\begin{equation}
  \langle  \eta, 1 \rangle=\int_{\Omega} \eta(x,y) dxdy =0,
\label{aver1}
\end{equation}
where $\langle u, v \rangle $ is the scalar product in $L_2(\Omega)$:
\begin{equation}
 \langle u, v \rangle =\int_0^h (\int_0^{2\pi} u(x,y) v(x,y) dx) dy.
\label{inprod}
\end{equation}
Note that if $u(x,y, t)$ is a solution to (\ref{OB2}),(\ref{boundNeum}) and (\ref{period2}),  then for  any constant $C$ the function
$u(x,y,t)  + C$ also is a solution.

We use below
the stream function - vorticity   formulation of these
equations in order to exclude the pressure $p$.
Introducing the vorticity $\omega$ and the stream function $\psi$, we obtain \cite{Chorin}
\begin{equation}
\Delta \psi=-\omega,
\label{OBEstream1}
\end{equation}
where the velocity  $\bf v$ can be expressed
 through the stream
function $\psi(x,y)$  by the relations
$v_1=\psi_y, v_2=-\psi_x$.
 Equations (\ref{OB1}), (\ref{div})
and (\ref{OB2})  take the form  \cite{Chorin}
\begin{equation}
\omega_t + \{\psi, \omega\} =\nu \Delta \omega,
\label{OBEstream2}
\end{equation}
here  $\{\psi, \omega \}=
  \psi_y \omega_x -\psi_x \omega_y$,
\begin{equation}
u_t + \{ \psi, u \} =\Delta u + \eta.
\label{heat1}
\end{equation}
The boundary conditions  become
\begin{equation}
  \psi(x,y,t)=\psi(x+2\pi,y,t),  \quad
\omega(x,y,t)=\omega(x+2\pi,y,t),
\label{boundstream1}
\end{equation}
\begin{equation}
 {\psi} (x,y,t)\vert_{y=h}   = \omega(x, y, t)\vert_{y=h}=0,
\label{boundstream2}
\end{equation}
\begin{equation}
 {\psi} (x,y,t)\vert_{y=0} =0, \quad   \omega(x, y, t)\vert_{y=0}= u_x(x,0,t).
\label{Maran1}
\end{equation}
\begin{equation}
  u_y(x, y, t)\vert_{y=h} =0, \quad u_y(x, y, t)\vert_{y=0} =0.
\label{Maran2}
\end{equation}
The vortex-stream reformulation of the problem is given by  equations and boundary conditions (\ref{OBEstream1}) - (\ref{Maran2}) and 
initial  conditions  
\begin{equation}
  u(x, y, 0) =u_0(x,y), \quad \omega(x, y, 0) =\omega_0(x,y).
\label{Maran3}
\end{equation}

\section{ Main result }
\label{sec:3}

Before to formulate  the main theorem, let us describe
 the method of the realization of vector fields (RVF)
invented by
 P. Pol\'a\v cik (see \cite{Pol1,Pol2}).
We change slightly the original
version to adapt it for our goals.

Let us consider a family of local semiflows $S^t_{\mathcal P}$ in a fixed Banach
space $B$. Assume
 these semiflows  depend on a parameter ${\mathcal P} \in B_1$, where
$B_1$  is another Banach space.
Denote by ${\mathcal B}^n$ the unit ball $\{q: |q|\le 1\}$ in ${\bf R}^n$, where $q=(q_1, q_2, ..., q_n)$ and
$ |q|^2=q_1^2 +
... + q^2_n$.  Remind that 
a  set $M$ is said to be locally invariant in an open set $W \subset B$  under a semiflow $S_t$ in $B$ if  $M$ is a subset of $W$ and each trajectories of $S_t$ leaving $M$ simultaneously leaves
$W$.
Consider a system of differential equations defined on the ball ${\mathcal B}^n$:
\begin{equation}
 \frac{dq}{dt}=Q(q),
\label{ordeq}
\end{equation}
 where
\begin{equation}
   Q \in C^1({\mathcal B}^n), \quad \sup_{q \in {\mathcal B}^n}|\nabla Q(q)| < 1.
\label{cond1}
\end{equation}
 Assume the vector field $Q$ is directed strictly
inward at the boundary $\partial {\mathcal B}^n=\{q: |q|=1 \}$:
\begin{equation}
   Q(q)\cdot q < 0 , \quad  q \in \partial {\mathcal B}^n.
\label{inward}
\end{equation}
Then system (\ref{ordeq}) defines a
global semiflow
 on ${\mathcal B}^n$. Let $\epsilon$ be a positive number.

\begin{definition} ({\bf realization of vector fields})  \label{RealVF} 
{We say that the family of local  semiflows $S^t_{\mathcal P}$  realizes the vector field $Q$ (dynamics (\ref{ordeq})) with accuracy $\epsilon$
(briefly, $\epsilon$  - realizes),
if there exists a  parameter ${\mathcal P}={\mathcal P}(Q, \epsilon, n)  \in B_1$
such that

({\bf i}) semiflow $S^t_{\mathcal P}$ has a positively invariant   and locally attracting
 manifold ${\mathcal M}_n \subset B$ diffeomorphic to the unit ball $ {\mathcal B}^n$;

({\bf ii})  this manifold is embedded into $B$  by a  map
\begin{equation}
  z = Z(q), \quad q \in {\mathcal B}^n, \quad z \in B, \quad Z \in C^{1+r}
({\mathcal B}^n),
\label{manifRVF}
\end{equation}
where $r > 0$;

({\bf iii}) the restriction of the semiflow  $S^t_{\mathcal P}$ to ${\mathcal M}_n$  is defined by the system of differential equations
\begin{equation}
  \frac{dq}{dt}=Q(q) + \tilde Q(q, {\mathcal P}), \quad Q \in C^{1}({\mathcal B}^n),
\label{reddynam}
\end{equation}
  where
\begin{equation}
	    |\tilde Q(\cdot, {\mathcal P})|_{C^1({\mathcal B}^n)} <
\epsilon.
\label{estRVF}
\end{equation}}
\end{definition}

This means that the dynamics on the invariant manifold
is defined by   the variables $q_1, q_2, ..., q_n$ and
approximates prescribed dynamics (\ref{ordeq}) with
accuracy $\epsilon$.
\vspace{0.2cm}

The  IBVP defined by  (\ref{OBEstream1}) -(\ref{Maran3}) involves the numbers
$ \nu,  h$ and the function $\eta(x,y)$.
We set
${\mathcal P}=\{h, \nu, \eta(\cdot,\cdot) \}$. The main result  is as follows:

\begin{theorem} \label{maint}
{ Dynamics of the semiflows defined by IBVP  (\ref{OBEstream1}) -(\ref{Maran3}) ,   is maximally complex in the following sense.
For each integer $n$, each $\epsilon > 0$  and each vector field
$Q$ satisfying (\ref{cond1}) and (\ref{inward}), there exists a value of the parameter
${\mathcal P}(Q,\epsilon)$ of this IBVP such that
this problem defines a  semiflow $S^t_{\mathcal P}$, which $\epsilon$-realizes
the vector field $Q$}.
\end{theorem}

This result implies the following corollary.

{\bf Corollary.} {\em  The family of semiflows  $S^t({\mathcal P})$ induced  by IBVP   (\ref{OBEstream1}) -(\ref{Maran3}) 
 with parameter $\mathcal P$, generate 
all  (up to  orbital topological equivalencies)   hyperbolic dynamics on compact invariant hyperbolic sets  defined by
$C^1$-smooth vector fields  on  finite dimensional smooth compact manifolds}.

\vspace{0.2cm}
In particular, we find that the Navier-Stokes dynamics can generate  Smale axiom A flows,
Ruelle-Takens attractors \cite{Ru2, NRT}, the Anosov flows, and, due to the persistence of compact invariant
hyperbolic sets \cite{Ru, Katok},  hyperbolic dynamics defined on these sets.

\section{ RVF method for weakly nonlinear evolution equations}
\label{sec:4}

This section describes a  general construction of the RVF method for "small" solutions.
Let us consider an  evolution equation
\begin{equation}
v_t=Lv + F(v) + \gamma f,
\label{eveq1}
\end{equation}
where $v$ lies in an Hilbert  space $H$, $L$ is a sectorial operator, $F$ is a nonlinear operator,
 $f \in H$ is independent of $v, t$ ("an external force") and  $\gamma > 0$ is a small parameter.
We use  the standard function spaces \cite{He}
$$
H_{\alpha}=\{v \in H:  ||v||_{\alpha}=||(I-L)^{\alpha} v|| < \infty\}.
$$
Assume  $F$ is a $C^{1+r}$ map from $H_{\alpha}$ to $H$, $r \in (0,1)$. We set
\begin{equation}
  v(0)=v_0, \quad v_0 \in H_{\alpha}.
\label{indata}
\end{equation}
We also suppose that   this map satisfies conditions
\begin{equation}
||F(v)|| \le  C_1||v||_{\alpha}^2,   \quad ||DF(v)|| \le  C_2||v||_{\alpha}
\label{FDF}
\end{equation}
for some $\alpha \in (0,1)$.

Then a unique solution of the Cauchy problem (\ref{eveq1}), (\ref{indata})  exists
on some open time interval $(0, t_0(v_0))$, $t_0 > 0$  \cite{He}.
The following  assumption plays a key role.
\vspace{0.2cm}

{\bf Spectral Gap Condition. }  {\em
Assume 
the spectrum $Spec \ L \subset {\bf C}$ of  $L$ consists of  the two parts: $Spec \ L=\{0\} \cup {\mathcal S}$,
where
\begin{equation}
  Re \ z  < - c_0 < 0 \quad for \ all \ z \in {\mathcal S}
\label{eveq3}
\end{equation}
and there exist exactly $N$ linearly independent $e_j \in H$ such that
\begin{equation} \label{eveq2}
  L  e_j=0, \quad  j=1,..., N. 
\end{equation}
}
\vspace{0.2cm}

Let $B_1$ be the space
$
B_1=Span \{e_1,..., e_N \}.
$
Then there exists a space $B_2$ invariant under $L$ such that $H=B_1 + B_2$, where $B_1+ B_2$  is a direct sum of $B_i$
(\cite{Danford1, Danford2}, also see \cite{He}, Th. 1.5.2).
We have two  complementary  projection operators ${\bf P}_1$ and ${\bf P}_2$ such that
   ${\bf P_1} +{\bf P_2}={\bf I}$, where
$\bf I$ denotes the identity operator, and $B_i={\bf P}_i H$. Let us denote by $L^*$ an operator  conjugate to $L$.
If the operator $L$ has a compact resolvent, then the spectra of $L$ and $L^*$ are discrete and we have  countable sets
of eigenvectors $e_j$  and $\tilde e_j$ of the operators $L$ and $L^*$ respectively, $j \in {\bf N}=\{1, 2, ... \}$.
In this case  without loss of generality one  can assume that $e_i$ and $\tilde e_i$ are biorthogonal:
$
\langle e_i, \tilde e_j \rangle =\delta_{ij},
$
where $\delta_{ij}$ stands for the Kronecker symbol.
Then $P_1$ can be defined by
\begin{equation}
 {\bf P_1} u =\sum_{i=1}^N \langle u, \tilde e_i \rangle e_i,
\label{PP1}
\end{equation}
where $
L^* \tilde e_i=0,  \quad i=1, ..., N.
$

Consider small  solutions of  (\ref{eveq1}) of the following form:
\begin{equation}
 v=\gamma X + w, \quad w(t) \in B_2,   \quad X(t)=\sum_{i=1}^N X_i(t) e_i   \in B_1.
\label{eveq4}
\end{equation}
Substituting (\ref{eveq4}) to eq. (\ref{eveq1}) we obtain
\begin{equation}
X_t= \gamma^{-1} {\bf P}_1 F(\gamma X + w) + f_1,
\label{eveq5}
\end{equation}
\begin{equation}
w_t=L w +  {\bf P}_2 F(\gamma X + w) + \gamma f_2,
\label{eveq6}
\end{equation}
where  $f_k={\bf P}_k f$.
Assume
$$
f_1=\gamma \bar f_1, \quad ||\bar f_1|| < \bar C_1, \quad || f_2|| < C_2.
$$
Let us set
$$
w=\gamma w_0 +\tilde w
$$
where  $w_0$ is defined by
$$
w_0 = - L^{-1} f_2.
$$
We consider (\ref{eveq5}), (\ref{eveq6}) in the domain
\begin{equation}
{\mathcal D}_{R_1, \gamma, C}=\{ (X, \tilde w): \ |X| < R_1, \ ||\tilde w||_{\alpha} < C \gamma^2 \}.
\label{domeveq}
\end{equation}
Let  ${\mathcal B}^N(R_0)$ be the ball ${\mathcal B}^N(R_0)=\{X: \ |X| < R_0 \}$.
The following   assertion will be useful below.

\begin{lemma} \label{rvf}
{Assume $r \in (0,1)$, $\alpha \in (0,1)$  and $C > C_0(R_1, F, \alpha)$ is large enough. Then
for sufficiently small positive $\gamma  < \gamma_0(r, \alpha, F, C,  R, N)$ system (\ref{eveq5}), (\ref{eveq6}) has
a locally invariant in ${\mathcal D}_{R_1, \gamma, C}$ and locally attracting  manifold ${\mathcal M}_{N, \gamma}$ defined by
\begin{equation}
w=\gamma (w_0 +  \tilde W(X, \gamma)),
\label{wW}
\end{equation}
where a $C^{1+r}$ smooth map $\tilde W: {\mathcal B}^N(R) \to H_{\alpha}$ satisfies the estimates
\begin{equation}
\sup_{X \in {\mathcal B}^N(R)} ||\tilde W(\cdot, \gamma)||_{\alpha}  +  \sup_{X \in {\mathcal B}^N(R)}  ||D_X \tilde W(\cdot, \gamma)||_{\alpha} <  c_1\gamma.
\label{tW1}
\end{equation}
}
\end{lemma}

{\bf Proof}.
This assertion is an immediate consequence of Theorem 6.1.7  \cite{He}. The proof is standard and can be found in  Appendix 2.

On the manifold ${\mathcal M}_{N, \gamma}$   evolution equation (\ref{eveq5}) for the slow component $X$ takes the following form:
\begin{equation}
\frac{dX}{d\tau}=  { K}(X,\gamma)  + { M(\gamma)} X  + \hat f_1 + \phi(X, \gamma),
\label{eveq8}
\end{equation}
where $\tau=\gamma t$,
\begin{equation}
{ K}(X, \gamma)={\bf P}_1 \gamma^{-2}(F (\gamma(X+w_0))  - \gamma DF(\gamma w_0)  X - F(\gamma w_0)),
\label{QX}
\end{equation}
 and ${ M(\gamma)}: {\bf R}^N \to {\bf R}^N$ is a bounded linear operator defined by
\begin{equation}
{ M}(\gamma) X = {\bf P}_1 \gamma^{-1} DF(\gamma w_0) X.
\label{MX}
\end{equation}
We  have
\begin{equation}
\hat f_1= \bar f_1  + {\bf P}_1 \gamma^{-2} F(\gamma w_0),
\label{hatf}
\end{equation}
and $\phi$ is a small correction such that
\begin{equation}
|\phi|_{C^{1+r}({ \mathcal B}^N(R))} < c_5 \gamma, \quad r >0.
\label{eveq91}
\end{equation}
For quadratic nonlinearities $F$  such that
$$
F(\alpha v)=\alpha^2 F(v)
$$
the relations for $ K$ and $ M$ can be simplified ( $K$ and $M$ do not depend on $\gamma$):
\begin{equation}
{ K}(X)={\bf P}_1 (F (X+w_0)  - DF( w_0) X - F( w_0)),
\label{QX1}
\end{equation}
\begin{equation}
{ M}(w_0) X = {\bf P}_1  DF(w_0) X.
\label{MX1}
\end{equation}

The key idea  is to consider the operator $M$ as a parameter in the RVF method.
This idea works if the following property holds.
\vspace{0.2cm}

{\bf Linear operator density (LOD) condition. }
{\em
Let us consider
 the set ${\bf O}_F$ of all linear operators  ${\bf M}(L^{-1} f)$ that can be obtained by (\ref{MX1}) when  $f$ runs
over the whole  space $B_2$.  We assume that this set  ${\bf O}_F$
is dense in the set of all linear operators ${\bf R}^N \to {\bf R}^N$.}
\vspace{0.2cm}

In coming sections  we apply this general approach to  IBVP   (\ref{OBEstream1})-(\ref{Maran3}).

\section{Existence and uniqueness for IBVP  (\ref{OBEstream1})-(\ref{Maran3}) and auxiliary estimates.}
\label{sec:7}

\subsection{ Function spaces and embeddings}
\label{sec: 7.1}

We  use standard Hilbert  spaces \cite{He}.
 We denote by $H=L_2(\Omega)$ the Hilbert space  of measurable, $2\pi$- periodical in $x$ functions defined
on $\Omega$
 with bounded  norms $|| \ ||$, where $||u||^2=\langle u, u\rangle$ and $\langle, \rangle$ is
the inner product defined by (\ref{inprod}).
Let us  denote by  $H_{\alpha}$ the fractional spaces
\begin{equation}
H_{\alpha}= \{ \omega :  ||\omega||_{\alpha} =||(I-\Delta_D)^{\alpha} \omega|| < \infty \},
\label{Hal}
\end{equation}
here $\Delta_D$ is the Laplace operator  with the standard domain corresponding to the zero Dirichlet boundary conditions:
\begin{equation}
Dom \ \Delta_D=\{ \omega: \omega \in W_{2,2}(\Omega), \quad \omega(x, y)\vert_{y=0, y=h}=0 \},
\label{DomDelta}
\end{equation}
here $W_{q,2}(\Omega)$ denote the standard Sobolev spaces.
 Let  $\tilde H_{\alpha}$  be another fractional space associated with $L_2(\Omega)$:
\begin{equation}
\tilde H_{\alpha}= \{ u :  ||u||_{\alpha} =||(I-\Delta_N)^{\alpha} u|| < \infty \},
\label{Bp}
\end{equation}
where $\Delta_N$ is the Laplace operator with the  domain corresponding to the zero Neumann boundary conditions
\begin{equation}
Dom \ \Delta_N=\{u: u \in W_{4,2}(\Omega), \quad u_y(x, y)\vert_{y=0, y=h}=0 \}.
\label{DomDeltaN}
\end{equation}
We  omit sometimes the indices $N, D$.  This choice of the domain for $\Delta_N$ is connected with
a special choice of main function space for $u$ -component.
The  Sobolev embeddings
\begin{equation}
H_{\alpha} \subset C^{s} (\Omega), \quad 0 \le s < 2(\alpha- 1/2),
\label{emb1}
\end{equation}
and
\begin{equation}
H_{\alpha} \subset L_q (\Omega), \quad 1/q > 1/2- \alpha, \ q \ge 2
\label{emb2}
\end{equation}
are useful below,
the same embeddings hold for $\tilde H_{\alpha}$.
Let $Tr(u)$ be the trace of a function $u$ on the bottom boundary $S$ of $\Omega$ (where $y=0$). We shall use  the following
embedding \cite{Stein}:
 \begin{equation}
|| Tr (u)||_{L_2(S)} \le c || u||_{\alpha}, \quad  \alpha > 1/4.
\label{Tr}
\end{equation}

\subsection{  Some preliminaries and auxiliary estimates}
\label{sec:7.2}

In coming subsections, our aim is to prove that IBVP  (\ref{OBEstream1})-(\ref{Maran3}),   defines a local semiflow. Moreover, we need
some estimates important for the invariant manifold technique.
To show existence of solutions we  use the standard semigroup methods.
Here, however, we met some difficulties because the Marangoni condition induces a singularity \cite{Pardo}.
To circumvent them, we choose  an appropriate phase space  taking into account that
the  temperature field $u$ should be more regular than the vorticity field $\omega$.  Moreover, we  use a special representation of the vorticity in order to represent
  IBVP (\ref{OBEstream1})-(\ref{Maran3})   by evolution equations.   

Let us consider IBVP  (\ref{OBEstream1})-(\ref{Maran3})  in the phase space
\begin{equation} \label{phasespace}
{\mathcal H}=H \times \tilde H_{1}, 
\end{equation}
 i.e., we use the space $\tilde H_{1}$ for $u$-component and the space $H$
for $\omega$ -component.
  
In order to apply semigroup technique,   we represent $ \omega$ as a sum of two terms, $\omega= \bar \omega +\tilde \omega$,
where the second term $\tilde \omega$ satisfies  the zero Dirichlet boundary conditions,
 and $\bar  \omega$ is defined as a solution of the following linear boundary value  problem:
\begin{equation}
   \Delta \bar \omega =0,
\label{Stokes}
\end{equation}
\begin{equation}
\bar \omega(x,0, t)= u_x(x,0,t), \quad \bar\omega(x, h, t)=0.
\label{omegab}
\end{equation}
We have $\psi=\bar \psi +\tilde \psi$, where
the functions $\bar \psi$ and $\tilde \psi$ are defined as solutions of the boundary value problems
 \begin{equation}
\Delta \bar \psi=- \bar \omega,
\label{psieq}
\end{equation}
 \begin{equation}
 \bar \psi(x, h, t)=\bar \psi(x, 0,t)=0,
\label{psib}
\end{equation}
\begin{equation}
\Delta \tilde \psi=- \tilde \omega,
\label{tpsieq}
\end{equation}
 \begin{equation}
 \tilde \psi(x, h, t)=\tilde \psi(x, 0, t)=0.
\label{tpsib}
\end{equation}
These problems can be resolved by the Fourier series.
It is  clear that, for sufficiently smooth $u$,  boundary value problem
(\ref{Stokes}),(\ref{omegab})
defines a linear operator $u \to \bar \omega(u)$. The following lemma gives useful estimates of this operator.

\begin{lemma} \label{5.1} {  The map  $u \to \bar \omega(u)$  satisfies
\begin{equation}
   ||\bar \omega(u)|| \le c_1|| u||_{\alpha}, \quad  \alpha > 1/2,
\label{Stokes1}
\end{equation}
\begin{equation}
   ||\nabla \bar \omega(u)|| \le ||\bar \omega||_{1/2} \le c_2|| u||_{\alpha_1}, \quad  \alpha_1 > 1,
\label{Stokes2}
\end{equation}
and for solutions $\bar \psi$ of (\ref{tpsieq}), (\ref{tpsib}) one has
\begin{equation}
   \sup |\nabla \bar \psi| \le c_3 ||u||_{\alpha}, \quad ||\nabla \bar \psi || \le c_3 ||u||_{\alpha}, \quad \alpha >3/2.
\label{Stokes3}
\end{equation}
}
\end{lemma}

We prove the lemma using  explicit solutions of boundary value problems (\ref{Stokes}),(\ref{omegab})
and (\ref{psieq}), (\ref{psib}) by the Fourier series.
The function $\bar\omega$ can be expressed via the Fourier coefficients of the trace $u$ on $S$ (see Appendix 1, (\ref{Ap4})). This gives
$$
||\bar \omega|| \le c_1 || Tr((I - D_x^2)^{1/4} u) ||_{L_2(S)}, \quad D_x=\frac{\partial }{\partial x}.
$$
By (\ref{Tr}) one obtains
(\ref{Stokes1}). Similarly,
$$
|| \nabla \bar \omega|| \le c || Tr((I - D_x^2)^{3/4}u) ||_{L_2(S)},
$$
and  embedding (\ref{Tr}) implies (\ref{Stokes2}). Estimate (\ref{Stokes3}) also follows
from the Fourier series (see Appendix 1).

\subsection{Evolution equations for $\tilde \omega$ and $u$}

Let us introduce some auxiliary maps.  
Let $\psi(\omega(\cdot, \cdot, \cdot))$ be a linear functional of $\omega$ defined as a periodic in $x$ solution of the boundary problem
\begin{equation}
 \Delta \psi=-\omega,    \quad \psi(x, 0,t)=\psi(x, h, t)=0.
\label{psi0}
\end{equation}
Below to simplify notation we denote  $\psi(\omega(\cdot, \cdot, \cdot))$ simply $\psi(\omega)$ or $\psi$.

The  map  $G: \omega \to G(\omega)$  is defined on $H$ by
\begin{equation} \label{mapG} 
  G(\omega)=-\{ \psi(\omega), \omega \}.
\end{equation}
The map  $F$  is defined on $H \times \tilde H_1$  by  
\begin{equation} \label{mapF}   
F(\omega, u)=-\{ \psi(\omega), u  \}.
\end{equation}

Moreover,  we introduce
\begin{equation}  \label{mapZ}
Z(\tilde \omega, u)=G(\tilde {\omega} +  \bar {\omega}(u)) -  \bar {\omega}(F).
\end{equation}
One has
$$
\bar \omega(u_t)=\bar \omega (\Delta u + F(\bar \omega +\tilde \omega, u)).
$$
Using this relation and definitions  (\ref{mapG}), (\ref{mapF}) and (\ref{mapZ}) we rewrite  IBVP defined by  (\ref{OBEstream1}) -(\ref{Maran3}) as a system of evolution equations:
\begin{equation} \label{equ}
 u_t=L_2 u  + F(\tilde \omega +  \bar \omega(u), u) +\eta,
\end{equation}
 \begin{equation} \label{eqw}
 \tilde {\omega}_t =L_1 \tilde \omega    + Z(\tilde \omega, u),  
\end{equation}
  where
$L_2=\Delta_N$ is the Laplace operator under the zero Neumann boundary conditions
(see (\ref{DomDeltaN})), and the linear operator $L_1$ is defined by
\begin{equation}
  L_1 \tilde \omega=\nu \Delta_D \tilde \omega - \bar \omega(\Delta_N u).
\label{L1}
\end{equation}
 Note that boundary conditions  (\ref{boundNeum}) and others  are taken into account by domain definitions
for operators $\Delta_D$ and $\Delta_N$.

\subsection{Smoothness of maps $F$ and $Z$}

We define
 spaces ${\mathcal H}_{\alpha}$  with $\alpha \ge 0$ and ${\mathcal H}={\mathcal H}_0$ by
\begin{equation}
  {\mathcal H}_{\alpha} = H_{\alpha} \times \tilde H_{\alpha + 1}, \quad {\mathcal H} = H \times \tilde H_{1}.
\label{Gest}
\end{equation}
In order to apply the  semigroup approach \cite{He} to system (\ref{equ}, \ref{eqw}), 
let us show first that $F: (\tilde \omega, w) \to F(\tilde \omega, w)$ and $Z: (\tilde \omega, w)\to Z(\tilde \omega, w) $ are $C^1$ smooth maps from  the space ${\mathcal H}_{\alpha}$ to the spaces $\tilde H_1$
and $H$, respectively.

\begin{lemma} \label{5.2} { $F$ defines a bounded map  from   ${\mathcal H}_{\alpha}$ to $\tilde H_1$}.
\end{lemma}

To prove this assertion, we use the estimate
\begin{equation}
  ||F(\tilde \omega, u)||_{1} \le S_1 + S_2 + S_3 + S_4,
\label{Fest1}
\end{equation}
where
$$
S_1=||\{(I-\Delta) \tilde \psi, u\}||, \quad S_2=2||\nabla \tilde \psi_x \cdot \nabla u_y - \nabla \tilde \psi_y \cdot \nabla u_x||,
$$
$$
S_3=||\{\tilde \psi, \Delta u\}||, \quad S_4=||\{\tilde \psi, u\}||.
$$
By  embeddings (\ref{emb1}) one has
$$
  S_1 \le c_1 ||\nabla  \Delta \tilde \psi)|| \sup |\nabla u| \le c_1||\tilde \omega||_{1/2} || u||_{\alpha+1}, \quad \alpha >0.
$$
In a similar way,
  $$
  S_2 \le c_2  (||\nabla \tilde \psi_x|| + ||\nabla \tilde \psi_y||) | u|_{C^2} \le c_3 ||\tilde \omega|| || u||_{\alpha+1}, \quad \alpha > 1/2.
$$
To estimate $S_3$ we use the inequalities
$$
S_3 \le c_4 ||\nabla \Delta u|| \sup |\nabla \tilde \psi| \le c_5 ||u||_{3/2} || \tilde \omega||_{\alpha},  \quad \alpha > 1/2.
$$
For $S_4$ one has
$$
S_4 \le  2||\nabla \tilde \psi|| \sup|\nabla u| \le c_6  ||\tilde \omega|| || u ||_{\alpha}.
$$
Therefore,
\begin{equation}
  ||F(\tilde \omega, u)||_{1}
 \le c ||\tilde \omega||_{\alpha} ||u ||_{1+\alpha}, \quad \alpha \in (1/2, 1).
\label{Fest3}
\end{equation}

Let us estimate $||F(\bar \omega, u)||_{1}$. Again one has
\begin{equation}
  ||F(\bar \omega, u)||_{1} =\bar S_1 + \bar S_2 + \bar S_3 + \bar S_4,
\label{Fest1b}
\end{equation}
where
$$
\bar S_1=||\{\bar \psi-\Delta \bar \psi, u\}||, \quad \bar S_2=2||\nabla \bar \psi_x \cdot \nabla u_y - \nabla \bar \psi_y \cdot \nabla u_x||,
$$
$$
\bar S_3=||\{\bar \psi, \Delta u\}||, \quad  \bar S_4=||\{\bar \psi, u\}||.
$$
Repeating  the same arguments as above and using Lemma \ref{5.1}, one has
$$
  \bar S_1 \le c_1 (||\nabla \Delta \bar \psi|| +||\nabla \bar \psi||) \sup |\nabla u| \le c_2||\bar \omega||_{1/2} || u||_{\alpha+1} \le
  c_3 || u||_{\alpha+1}^2, \quad \alpha >0.
$$
In a similar way,  one obtains
  $$
  \bar S_2 \le c_4  (||\nabla \bar \psi_x|| + ||\nabla \bar \psi_y||) | u|_{C^2} \le c_5 || \bar \omega|| || u||_{\alpha+1}
  \le c_6 || u||_{\alpha+1}^2 , \quad \alpha > 1/2.
$$
Moreover,
$$
\bar S_3 \le c_7 ||\nabla \Delta u|| \sup |\nabla \bar \psi| \le c_7 ||u||_{3/2} || \bar \omega||_{\alpha}
\le c_8 ||u||_{\alpha+1}^2,  \quad \alpha > 1/2
$$
and
$$
\bar S_4 \le  2||\nabla \bar \psi|| \sup| \nabla u|  \le c_9 ||\bar\omega|| ||u||_{\alpha+1} \le c_{10} ||u||_{\alpha+1},
\quad \alpha > 1/2.
$$
Therefore,
\begin{equation}
  ||F(\bar \omega, u)||_{1} \le c_{7} ||u||_{\alpha+1}^2, \quad \alpha > 1/2.
\label{Fest4}
\end{equation}
Combining this estimate with (\ref{Fest3}) one concludes that $F$ is a bounded
 map from ${\mathcal H}_{\alpha} $ to $\tilde H_1$.

\begin{lemma} \label{5.3} { $Z$ defines a bounded map  from   ${\mathcal H}_{\alpha}$ to $ H$}.
\end{lemma}

To prove this lemma, let us notice
\begin{equation}
  ||\{\psi,\omega\}|| \le c_1 (\sup|\nabla \tilde \psi| + \sup|\nabla \bar \psi|) (||\tilde\omega||_{1/2} + ||\bar \omega||_{1/2}).
\label{Gest1}
\end{equation}
This gives, by (\ref{Stokes3}), Lemma \ref{5.1},
and the Sobolev embeddings (\ref{emb1})  that
\begin{equation}
  ||\{ \psi,\omega\}||  \le c_4(||u||_{\alpha+ 1}^2+ ||\tilde \omega||_{\alpha}^2), \quad
  \alpha  > 1/2.
\label{Gest2}
\end{equation}
Let us   estimate  $\bar \omega(F)$. By Lemma \ref{5.1} one has
\begin{equation}
  ||\bar \omega(F)|| \le c_1|| F||_{\gamma}, \quad  \gamma > 1/2.
\label{Gest11}
\end{equation}
One has   $|| F||_{\gamma} \le c_2||F||_1$. The  estimate of $||F||_1$  follows from  Lemma \ref{5.2}. Therefore,
\begin{equation}
  ||\bar \omega(F)|| \le c_3(|| u||_{\alpha + 1}^2 + ||\tilde \omega||_{\alpha}^2), \quad  \alpha  > 1/2.
\label{Gest6b}
\end{equation}
The proof is complete.

Lemmas \ref{5.2} and \ref{5.3} show that
 $F$ and $Z$ are  bounded maps from a bounded domain in ${\mathcal H}_{\alpha}$ to
$\mathcal H$. They are quadratic and analogous estimates imply that the derivatives of $F, G$
are bounded as maps from
${\mathcal H}_{\alpha}$ to ${\mathcal H}$. Indeed, the derivative of $D_{\tilde \omega} G$ with respect to $\tilde \omega$
is a linear operator from $H_{\alpha}$ to $H$
defined by
$$
(D_{\omega} G)\delta \tilde \omega=  \{ \delta \tilde \psi, \omega \} + \{\psi, \delta \tilde \omega\}.
 $$
 An estimate of the norm of this operator can be found as above.

\subsection{ Transformation of evolution equations (\ref{equ}), (\ref{eqw}) and associated linear operator}
\label{sec: 7.4}

We follow the standard approach developed for the Rayleigh- B\'enard and Marangoni -B\'enard convection
\cite{Drazin,GJ, Bracard}. Assume that
the temperature field $u$ is a small $\gamma$ -perturbation of a vertical profile $U(y)$. Here  $\gamma >0 $ is a small parameter
independent of the viscosity $\nu$ (this assumption is  important): $\gamma < \gamma_0(\nu)$.
Let 
$U$ be a $C^{\infty}$-smooth function of $y \in [0,h]$ such that 
\begin{equation}  \label{suppU}
U(y)=U_y(y) =0,  \quad \forall y \in [0,\delta_1),
\end{equation}
for some $\delta_1  \in (0, h)$.
Assume   $u_1(x,y)$ is
a $C^{\infty}$ smooth  $2\pi$ periodical in $x$ function  satisfying conditions
\begin{equation}
 u_1(x,y)=0  \quad \forall x \in (-\infty, +\infty), \  \forall y \in (\delta_1, h].
\label{suppu1}
\end{equation}
We set
$$
 u_0 =U + \gamma u_1
$$
and
$$
\eta=\eta_0 + \gamma^2 \eta_1, \quad \eta_0=-\Delta u_0,
$$
$$
\langle \eta_1, 1 \rangle=0, \quad |\eta_1|_{C^3(\Omega)} < C_0.
$$
Let us represent $u$ as
\begin{equation}
u=u_0 + w,
\label{transform1}
\end{equation}
where $w$ a new unknown function.

For new unknowns $\tilde \omega$, $w$ system (\ref{equ}), (\ref{eqw}) takes the form
\begin{equation}
 \tilde \omega_t=\nu \Delta \tilde \omega  - \bar \omega(w_t) - \{\psi, \bar \omega + \tilde \omega \},
\label{eveq10}
\end{equation}
\begin{equation}
w_t= \Delta w   - \{\psi, U +\gamma u_1 + w \} + \gamma^2 \eta_1, 
\label{eveq11}
\end{equation}
where 
\begin{equation}
\omega=\bar \omega + \tilde \omega
\label{transform2}
\end{equation}
where the function $\bar \omega(x,y)$ is  the solution of boundary value problem (BVP)
(\ref{Stokes}), (\ref{omegab}) with $u=w$  (we can set $u=w$  due to above assumptions
on the supports of $u_1$ and $U$). Thus this BVP has the form
\begin{equation}
\Delta \bar \omega=0, \quad \bar \omega(x,h, t)=0,  \quad \bar \omega(x,0,t )=w_x.
\label{eveq13}
\end{equation}
This boundary value problem (BVP)  defines a linear map $w \to \bar \omega(w)$.

Removing the nonlinear terms and  ones of  order $\gamma$ and $\gamma^2$ in (\ref{eveq10}), (\ref{eveq11}),
we obtain the linear operator
\begin{equation}
 Lv=(\bar L_1 v, \bar L_2 v)^{tr},   \quad v=(\tilde \omega, w)^{tr}
\label{Lop}
\end{equation}
 where the  operators $\bar L_k$ are defined by
\begin{equation}
\bar L_1 v=\nu \Delta \tilde \omega  - \bar \omega(\bar L_2 v),
\label {A1}
\end{equation}
\begin{equation}
\bar L_2 v=\Delta w +  \psi_x U_y.
\label {A2}
\end{equation}
Here $\psi(\tilde \omega, w)$ is defined  by (\ref{psi0}) with $\omega=\bar\omega(w) + \tilde \omega$ in the right hand side.
Below we use  a natural decomposition
$\psi=\bar \psi + \tilde \psi$, where $\bar \psi, \tilde \psi$ are defined by
(\ref{psieq}), (\ref{psib}) and (\ref{tpsieq}), (\ref{tpsib}), respectively.

The spectral problem for the operator $L$
\begin{equation}
\lambda \tilde \omega= \bar L_1 v, \quad \lambda w= \bar L_2 v
\label {A2S}
\end{equation}
can be represented in the  standard form.  Indeed,
the equation for the first component $v_1=\tilde \omega$ takes the form
$$
\lambda \tilde \omega=\nu \Delta \tilde \omega - \bar \omega(\bar L_2 v),
$$
and, since $\Delta \bar \omega=0$,   $\bar \omega(\bar L_2 v)=\lambda \bar \omega(w)$,
 spectral problem ({\ref{A2S}) becomes
\begin{equation}
\lambda \omega=\nu \Delta \omega,
\label {Sp1}
\end{equation}
\begin{equation}
\lambda w=\Delta w +  \psi_x U_y,
\label {Sp2}
\end{equation}
where the functions $w$ and $\omega$ satisfy the boundary conditions
\begin{equation}
w_y(x, y)\vert_{y=0, h}=0,  \quad   \omega(x,h)=0,  \quad  \omega(x,0)=w_x.
\label {Sp3}
\end{equation}
   This spectral problem  is investigated in coming sections but first we consider some general properties of $L$.

\subsection{ $L$ is a sectorial operator}
\label{sec: 7.5}

In order to apply the standard technique \cite{He}, first let us  show  that the operator $-L$
is sectorial. We use the following known result \cite{He}: if $L^{(0)}$ is a self adjoint operator in a Banach space $X$,
$L^{(0)}: X \to X$ and
$B$ is a linear operator, $B: X \to X$ such that $Dom \ L^{(0)} \subset Dom \ B$ and for all $\rho \in Dom \ L^{(0)}$
\begin{equation}
||B \rho || \le \sigma || L^{(0)} \rho || + K(\sigma) ||\rho||
\label{sect}
\end{equation}
for  $0 < \sigma < 1$ and some constant $K(\sigma)$, then $L^{(0)} + B$ also is a sectorial operator.

  Let us define the unperturbed operator $L^{(0)}$ by the relations
$$
  \bar L_1^{(0)} (\omega,w)^{tr}=  \nu \Delta \omega,  \quad \bar L_2^{(0)} (\omega, w)^{tr} = \Delta w
$$
with  domains defined by (\ref{DomDeltaN}) and (\ref{DomDelta}),  respectively.
The operator $B$ is given then by
$$
 B (\omega, w)^{tr} =(-\bar \omega(\bar L_2 v),  \psi_x  U_y)^{tr}.
$$
The spectral problem for
$L^{(0)}$ is given by (\ref{Sp1}), (\ref{Sp2}) and (\ref{Sp3}) with $U=0$.  The corresponding
 eigenfunctions are
$$
   w_{k,m}(x,y)=\cos(m\pi y h^{-1})\exp(ikx), \quad \omega_m=0, \quad k,m=0,1,..., 
$$
with the  eigenvalues  $\lambda_{k,m}=-m^2 \pi^2 h^{-2}-k^2$, and
$$
   w_{k,n}(x,y)=0,    \quad \omega_{k,n}=\sin(n\pi y h^{-1})\exp(ikx), \quad n=1,2,..., \quad  k=0,1,..., 
$$
with the  eigenvalues  $\bar \lambda_{k,n}=- \nu n^2 \pi^2 h^{-2}-k^2$.

\begin{lemma} \label{5.4} { Under condition (\ref{suppU}) $L$ is a sectorial operator}.
\end{lemma}

The  operator $L^{(0)}$
is self-adjoint, the spectrum is discrete
 and lies in the interval $(-\infty, 0)$
  Therefore,  $L^{(0)}$ is a sectorial.

Let us check estimate (\ref{sect}). First we estimate $\bar \omega(\bar L_2 v)$.
Since $U$ satisifies (\ref{suppU}), 
one has $\bar \omega(\bar L_2 v)=\bar \omega(\Delta w)$. Thus,  for $\alpha \in (1/2, 1)$ and for any $\sigma >0$ by Lemma \ref{5.1}
we obtain
\begin{equation}
||\bar \omega(\Delta w)|| \le c ||\Delta w||_{\alpha} \le c ||w||_{1+\alpha}  \le K(\sigma) ||w||_1 + \sigma||\Delta w||_1.
\label{sect5}
\end{equation}
Let us consider the second component $\psi_x U_y$ of $B$, where $\psi_x =\bar \psi_x + \tilde \psi_x U_y$. Since
$U(y)$ is a smooth function,  and $\tilde \psi$ is a solution of boundary value problem (\ref{tpsieq}), (\ref{tpsib}) one obtains
\begin{equation} \label{tpsiest1}
||\tilde \psi_x U_y||_1 \le c ||\tilde \omega_x|| \le \sigma||\Delta \tilde \omega|| + K_1(\sigma)||\tilde \omega||
\end{equation}
for some $K_1(\sigma)$.
The function $\bar \psi$ is a solution of BVP (\ref{psieq}), (\ref{psib}). Therefore, 
we can represent  $\bar \psi$ as
$$
\bar \psi(x,y)= \int_0^{2\pi} \Gamma(x-s, y) w_s(s, 0) ds,
$$
where $\Gamma$ is analytic   in $x-s, y$ for $y > \delta_1 >0$.  Since $U$ satisfies (\ref{suppU}) this relation entails 
$$
||\bar \psi_x U_y||_1 \le C_2 \int_0^{2\pi} |w(s,0)|^2 ds. 
$$
By estimate (\ref{Tr}) for traces  one has
\begin{equation} \label{compact}
||\bar \psi_x U_y||_1 \le C_3  ||w||_{1/2} \le \sigma ||w||_1 + K_2(\sigma)||w||. 
\end{equation}
This estimate together with  (\ref{tpsiest1}) and (\ref{sect5}),
 proves the lemma.

\subsection{Existence and uniqueness}

The fact that the linear operator $L$ is sectorial and
the $C^1$ smoothness of the maps $F$ and $Z$ entail  that
system (\ref{equ}) -(\ref{eqw}) defines a local $c^1$-smooth semiflow.
For results on global existence see \cite{Glob}. We do not use these results in this paper. 
The global existence on time interval $[0, +\infty)$ will be proved only for trajectories, which lie in a small open neighborhood of the positively invariant manifold
${\mathcal M}_n$ (see definition  \ref{RealVF}).

\subsection{Resolvent of $L$ is a compact operator}

Let us prove that the resolvent of $L$ is a compact operator.

\begin{lemma} \label{resolv} {For sufficiently large $\nu>0$ the resolvent $(L- \lambda)^{-1}$ is a compact operator from $\mathcal H $ to $\mathcal H$ for some $\lambda$.
}
\end{lemma}

{\bf Proof}.
Let us take $\lambda=\nu^2$. 
Consider the boundary value problem that defines the resolvent:
\begin{equation}
\lambda \omega-\nu \Delta \omega =  f,     \quad \omega(x,h)=0, \quad \omega(x,0)=w_x(x,0),
\label{res1}
\end{equation}
\begin{equation}
\lambda w-  \Delta w =  \psi_x U_y + g,     \quad w_y(x,y)\vert_{y=0, h}=0,
\label{res2}
\end{equation}
where
\begin{equation}
\Delta \psi=-\omega,  \quad  \psi(x,0)=0, \quad \psi(x, h)=0
\label{res3}
\end{equation}
and $f \in H$, $g \in H_1$, i.e.,
\begin{equation}
|| f|| + ||g||_1 < C_0.
\label{res13}
\end{equation}
We assume that
$
  \langle g, 1 \rangle =0$.
We represent $\omega$ as a sum $\omega=\bar \omega + \tilde \omega$, where
\begin{equation}
\lambda \tilde\omega- \nu \Delta \tilde\omega =  f,     \quad \tilde\omega(x,h)=0, \quad \tilde\omega(x,0)=0,
\label{res4}
\end{equation}
\begin{equation}
\lambda \bar \omega - \nu \Delta \bar \omega=0,      \quad \bar \omega(x,h)=0, \quad \bar \omega(x,0)=w_x(x,0).
\label{res5}
\end{equation}
In order to prove the lemma, it is sufficient to obtain for some $\alpha  >0$ the following estimates:
\begin{equation} \label{tomega}
||\tilde \omega||_{\alpha} \le c || f||,
\end{equation}
\begin{equation} \label{wfg}
||w ||_{1+ \alpha} \le  c_2 (|| f|| + ||g||_1),
\end{equation}
\begin{equation} \label{bomega}
||\bar \omega ||_{\alpha} \le c_1(|| f|| + ||g||_1).
\end{equation}

Since $\lambda=\nu^2$ the solution $\tilde \omega$ of eq. (\ref{res4})  satisfies estimate  (\ref{tomega}).
 Note that $\tilde \omega$ does not depend on $w$.  We use the natural decomposition 
$\psi=\tilde \psi+\bar \psi$ defined by (\ref{psieq}),(\ref{psib}),    (\ref{tpsieq}) and (\ref{tpsib}).  
Then BVP (\ref{res2}) can be rewritten in an operator form  as
\begin{equation}
w=A_{\nu}w + \tilde g,     
\label{res2a}
\end{equation}
where 
$$
A_{\nu}w=(\nu^2 - \Delta)^{-1} (\bar \psi_x(w(\cdot,\cdot)) U_y)  
$$
defines  a linear operator $A_{\nu}: H_1 \to H_1$
and
$$
\tilde g=(\nu^2 - \Delta)^{-1}  (\tilde \psi_x U_y +  g).
$$

Let us prove that $A_{\nu}$ is a contraction in the space $H_1$ for sufficiently large $\nu$.  
We use  estimate (\ref{compact}) that gives
\begin{equation} \label{compact2}
||A_{\nu} w||_1  \le \nu^{-1}(\sigma ||w||_1 + K||w||) 
\end{equation}
for some  $K >0$, which does not depend on $\nu$. 
This estimate   shows that  the operator $A_{\nu}$ is a contraction in $H_1$ for sufficiently large $\nu$. Therefore, the solution of 
(\ref{res2a}) exists and satisfies  
\begin{equation} \label{compact3}
||w||_{1+\alpha} \le C||\tilde g||_{\alpha}.
\end{equation}
One has
\begin{equation} \label{compact4}
||\tilde g||_{\alpha}  \le c( ||\tilde \psi_x  U_y||_{\alpha} + ||g||_1).
\end{equation}
Note that $U_y$ is a smooth bounded function and the solution $\tilde \psi$ of  (\ref{res4}) satisfies $||\tilde \psi_x||_{\alpha} \le c_2||f||$ for some $c_2 >0$.
Therefore, 
\begin{equation} \label{compact5}
||\tilde g||_{\alpha}  \le c_3( ||f|| + ||g||_1).
\end{equation}
Estimates (\ref{compact3}), (\ref{compact4}) and (\ref{compact5}) prove   (\ref{wfg}).  Estimate (\ref{bomega}) follows from   (\ref{wfg}) by the embedding for traces.
 The lemma is proved.

This lemma implies, according to the well known result (see \cite{Kato}, Ch. III, Theorem 6.29),  that the spectrum of $L$ is discrete (it consists of isolated eigenvalues), each
eigenvalue has a finite multiplicity $n(\lambda)$,
and the resolvent ${\mathcal R}(\lambda)$ is a compact operator for all $\lambda$, where ${\mathcal R}(\lambda)$ is defined.  We investigate the spectrum in the next section.

\section{Spectrum of the main  linear operator} \label{Spectrum}

\subsection{ Some preliminaries} \label{sec:8.1}

Let us consider  spectral problem (\ref{Sp1}), (\ref{Sp2}) and (\ref{Sp3}).  
For any $U(y)$ this problem has the trivial eigenfunction  $e_0=(0,1)^{tr}$, where $\omega=0, w=1$,  with
the zero eigenvalue $\lambda$.
We consider eigenfunctions $e(x,y,\lambda)$ with eigenvalues $\lambda \in {\bf  C}_{1/2}$, where 
 ${\bf C}_a$ denotes the half-plane 
\begin{equation}
{\bf C}_a= \{ \lambda \in {\bf C}:   Re \ \lambda > -a \}.
\label{Ck}
\end{equation}
Since $U$  depends only on $y$, 
we  seek the eigenfunctions  of the form 
\begin{equation}
   \quad \quad w(x,y, \lambda)= w_k(y, \lambda) \exp(ikx),
\label{wF}
\end{equation}
  \begin{equation}
   \psi(x,y, \lambda)= \psi_k(y, \lambda)  \exp(ikx), \quad   \omega(x,y, \lambda)= \omega_k(y,\lambda)  \exp(ikx).
\label{psiF}
\end{equation}
For  $\omega_k$, $\psi_k$ and $w_k$ one obtains the following
boundary value problem:
\begin{equation} \label {BVK1}
\frac{\partial^2 \omega_{k}}{\partial y^2} - k_{\nu}^2 \omega_k=0,   \quad  \omega_k(h,\lambda)=0, \quad \omega_k(0, \lambda)= ik w_k(0, \lambda),
\end{equation}
where $k_{\nu}^2=k^2 +\lambda/\nu$,
\begin{equation} \label {BVK2}
\frac{\partial^2 \psi_{k}}{\partial y^2} - k^2 \psi_k=-\omega_k,   \quad  \psi_k(h, \lambda)=0, \quad \psi_k(0, \lambda)=0,
\end{equation}
\begin{equation} \label {BVK3}
\frac{\partial^2 w_{k}}{\partial y^2}  - \bar k^2 w_k= ik U_y(y) \psi_k,   \quad  \frac{\partial w_k(y, \lambda)}{\partial y}\vert_{y=0, h}=0,
\end{equation}
where $\bar k=\sqrt{k^2 +\lambda}$.
Let us suppose, without loss of generality, that $k > 0$, and $Re \  \bar k >0$ for $\lambda\in {\bf C}_{1/2}$,   since $w_{-k}$ are functions, complex conjugate to $w_k$
and $\bar k$ is involved only via $\bar k^2$.  We assume that
\begin{equation} \label{hviscos}
h=10 \log \nu, \quad  \nu >> 1.
\end{equation}

The solution of  problem (\ref{BVK1}),  (\ref{BVK2}) is defined
by
 \begin{equation}
   \omega_k(y,\lambda) = \beta_k \frac{\sinh (\bar k_{\nu}(h-y))}{\sinh(\bar k_{\nu} h)},    \quad  \beta_k(\lambda)= ik w_k(0, \lambda) ,
\label{omFk}
\end{equation}
\begin{equation}
 \psi_k(y, \lambda)=- \nu \beta_k \lambda^{-1} \Phi_k(y, \lambda),
\label{psiF2}
\end{equation}
\begin{equation}
 \Phi_k(y, \lambda)= \frac{\sinh(k h) \sinh(\bar k_{\nu}(h-y)) -  \sinh(\bar k_{\nu}h)  \sinh (k(h-y)) }{ \sinh(\bar k_{\nu} h) \sinh(kh)}.
\label{psiF2a}
\end{equation}

Note that relation (\ref{psiF2}) is correctly defined for all $\lambda \in {\bf C}_{1/2}$,  in particular, for  $\lambda=0$. Indeed, ,  for     
small $\lambda$ 
\begin{equation}
\bar k_{\nu}-k=\sqrt{k^2 + \lambda \nu^{-1}} -k=  \lambda (2\nu k)^{-1} +  O(\lambda^2\nu^{-2} k^{-3})
\label{dk}
\end{equation}
that gives 
\begin{equation}
   \psi_k(y, \lambda)=\beta_k(\lambda) \frac{ y\sinh(kh) \cosh(k(h-y)) - h \sinh(ky) }{2k\sinh^2( k h)}  + \phi(y, k),
\label{psiF4}
\end{equation}
where
$$
|\phi(y, k, \lambda)| <  c|w_k(0, \lambda)||\lambda| (k\nu)^{-1}, \quad 0 < y < h.
$$
 For large $\nu_0$  and $|\lambda| << \nu$ assumptions  (\ref{hviscos})  allow us to simplify (\ref{omFk}) and (\ref{psiF2}). By (\ref{omFk})
 we obtain then  
\begin{equation} \label{omegaF}
  \omega_k(y, \lambda)= \beta_k(\lambda)  (\exp(- k y)  +  \tilde \omega_k(y, \nu)),
\end{equation}
where for each $s \in (0,1)$ and $|\lambda| < \nu^s$
\begin{equation} \label{tildeomega}
|\tilde \omega_k(y, \nu)| < C_s (|\lambda| (k \nu)^{-1} \exp(- ky) +  \exp(-k h)), \quad y \in [0,h],  
\end{equation}
where $C_s>0$ are constants independent of $s$ and $k$.
This estimate and (\ref{psiF4}) give
\begin{equation} 
\label{psiFa}
\psi_{k}(y, \lambda)=\beta_k(\lambda)(\bar \psi_k(y, \lambda) + \tilde \xi_k(y,  \lambda)),
\end{equation}
where 
\begin{equation} 
\label{psiFa1}
\bar \psi_k(y, \lambda)=\frac{y}{2k} \exp(- k y) ,
\end{equation}
and for $|\lambda| < \nu^s$
\begin{equation} \label{tildepsi}
|\tilde \xi_k(y,  \lambda)| <  \bar C_s (|\lambda| k^{-1} \nu^{-1} \exp(-ky) +    \exp(-k h)),   \quad y \in (0,h),
\end{equation}
 where constants $\bar C_s>0$ are uniform in $k, \nu$.

To investigate (\ref{BVK3}), we apply a  lemma.

\begin{lemma} \label{6.1}
{ Let us consider the  boundary value problem on $[0, h]$ defined by
\begin{equation}
w_{yy} - \bar k^2 w= f(y), \quad y \in [0, h],
\label{kwd1}
\end{equation}
\begin{equation}
w_{y}(y)\vert_{y=0, h} =0.
\label{kwb}
\end{equation}
Then
\begin{equation}
w(0)=-\int_0^h f(y)  \rho_{\bar k}(y) dy,
\label{kwc}
\end{equation}
where
\begin{equation}
  \rho_{\bar k}(y) = \frac{\cosh(\bar k(h- y))}{{\bar k} \sinh {\bar k} h}.
\label{kwr}
\end{equation}
}
\end{lemma}

To prove it, we multiply both the right hand and the left hand sides
of eq. (\ref{kwd1}) by $\rho_{\bar k}$ and  integrate by parts
in the left hand side. Note that  
\begin{equation} \label{rhoest1}
|\rho_{\bar k}(y)- \bar \rho_{\bar k}(y)| < \bar k^{-1} \exp(- \bar k h), \quad  \bar \rho_{\bar k}(y) =\bar k^{-1}\exp(-\bar k y). 
\end{equation}

\subsection{ Main result on spectrum of operator $L$}
\label{sec: 8.2}

Let us formulate the  assertion. 

\begin{proposition} \label{6.2}
{ Let  assumptions (\ref{hviscos}) hold,
 $N$ be a positive integer  and $K_N=\{1, ..., N \}  \subset {\bf Z}_+$.
Then 
there exists a  $C^{\infty}$ smooth  function $U(y)=U_{N}(y, \nu)$ satisfying (\ref{suppU}) and such that  for sufficiently large $\nu > \nu_0(N) >0$
the eigenfunctions $\lambda(k,  \nu)$ of BVP (\ref{BVK1}), (\ref{BVK2}) and (\ref{BVK3}) satisfy

({\bf i} )  
\begin{equation}
\lambda(k, \nu)=0  \quad    k \in K_N, 
\label{Spec0}
\end{equation}

({\bf ii}) 
\begin{equation}
  Re \ \lambda(k,  \nu) < - \delta_N  \quad  k \notin K_N,  
  \label{Spec}
\end{equation}
where positive $\delta_N$ is uniform in $\nu$.
}
\end{proposition}

{\bf  Proof}.
We use  Lemma \ref{6.1} to obtain a nonlinear equation for the eigenvalues $\lambda(k)$ of the boundary value
problem (\ref{BVK1})-(\ref{BVK3}). As a result,    one has
\begin{equation}
-k^2 \int_0^h   \psi_k(y,\lambda) \rho_{\bar k}(y) U_y(y, \nu)   dy=\beta_k(\lambda).
\label{kwc2}
\end{equation}

Note that the operator $L$ is not self-adjoint. Therefore, there are possible complex eigenvalues $\lambda$, i.e., complex roots of (\ref{kwc2}).
Moreover, let us note that, according to (\ref{psiF2}),  if $\beta_k( \lambda)=0$, then eq. (\ref{kwc2}) is satisfied.
In this case  (\ref{BVK1}) entails that
\begin{equation} \label {BVK10}
\frac{d^2\omega_{k}(y)}{dy^2} - k_{\nu}^2 \omega_k(y)=0,   \quad  \omega_k(h)=0, \quad \omega_k(0)=0,
\end{equation}
therefore $k_{\nu}^2=- (n \pi/h)^2$, where $n$ is an integer. This gives
$\lambda=-\nu((n \pi/h)^2+ k^2)< -1/2$. These eigenvalues $\lambda$ correspond to trivial solutions of the eigenfunction problem   with  $\lambda \notin {\bf C}_{1/2}$.  
Therefore,   without loss of generality we can set
$\beta_k(\lambda)=1$ in eq. (\ref{kwc2}).

The plan of the proof is as follows.
We consider the two  cases:
$
( {\bf I}) \hspace{0.1cm} |\lambda| < \nu^{3/4}
$
and
$
 ({\bf II}) \hspace{0.1cm} |\lambda| > \nu^{3/4}.
$
In the  first case we can simplify equation (\ref{kwc2}),
in the second case a rough estimate shows that
 eq. (\ref{kwc2})  has no solutions.

Let us start with the case {\bf I}.  To simplify our statement, 
we  first consider  a  formal limit of equation (\ref{kwc2}) as $\nu \to +\infty$.  Using ( \ref{psiFa1}) and (\ref{rhoest1}) one obtains that this limit has the form 
\begin{equation}
 \int_0^{+\infty}   y \bar U_y(y) \exp(-(k +\bar k) y)dy=2\bar k/k,
\label{kwc3}
\end{equation}
where  $\bar U=\lim  U_{K_N}(y, \nu)$ as $\nu \to +\infty$.
We  set $\bar U=V(y, d)$, where $V$ is a function  of a special form 
that depends on some parameters $d=(d_1,..., d_N)$.
Consider $C^{\infty}$ -mollifiers
$\delta_{\epsilon}(y)$ such that $\delta_{\epsilon} \ge 0$,  the support  $supp \ \delta_{\epsilon}(y)$ is
$(-\epsilon, \epsilon)$ and
\begin{equation}
\int_{-\epsilon}^{\epsilon}\delta_{\epsilon}(y)dy=1,
\label{Inteps}
\end{equation}
\begin{equation}
\sup |D_y^k\delta_{\epsilon}(y)| < c_k \epsilon^{-(k+1)}, \quad  k=0,1,2.
\label{Inteps1}
\end{equation}

Let us define  the function $V(y,d)$ on $(0, \infty)$ by 
\begin{equation}
 V_y(y, d)=2y^{-1}(\delta_{\kappa}(y-z_0) + \mu \chi(y-z_0) y  W_N(y, d)),   \quad  V(0, d)=0,
\label{Uy}
\end{equation}
where $\chi(z)$ is the step function  such that 
$
\chi(z)=1 \quad z >0$ and $\chi(z)=0 \quad z \le 0$,
$W_N$ is a polynomial in $y$ of the degree $N+1$  with coefficients depending on some parameters $d=(d_1, d_2, ..., d_N)$,
and
\begin{equation}
  \mu = \kappa^{2/3}, \quad z_0=5\kappa,     
\label{z0}
\end{equation}
where $\kappa$ is  a small parameter independent of $\nu$ as $\nu \to +\infty$. 
Coefficients of the polynomial $W_N$ and $d_j$ are assumed to be  bounded:
\begin{equation}
W_N(y, d)=\sum_{j=0}^{N} b_j(d) y^j, \quad  |b_j(d)| < C_j,   \quad  |d_j| < 1/2.
\label{PN}
\end{equation}
We set $U(y, \kappa)=\int_0^y V_s ds$ then $U$ satisfies condition (\ref{suppU}). 

To investigate (\ref{kwc3})
 it is useful to  introduce the variable
\begin{equation}
p=k +\bar k=k +(k^2 +\lambda)^{1/2}.
\label{pk}
\end{equation}
 Then  eq. (\ref{kwc3}) can be rewritten as
 \begin{equation}
   \frac{p}{k}=2  + S(p, k, d),
\label{ME}
 \end{equation}
 where
\begin{equation}
    S(p, k, d)=\mu G(p, d)  + g_{\kappa}(p)  + \tilde g_{\kappa}(p) ,
\label{kwe70}
 \end{equation}
 \begin{equation}
  g_{\kappa}(p)=  -1  + \int_0^{\infty} \delta_{\kappa}(y-z_0) \exp(-py) dy,
\label{kweg}
 \end{equation}
\begin{equation}
  \tilde g_{\kappa}(p)=  - \mu \int_0^{z_0} y W_N(y,d) \exp(-py) dy,
\label{kwegt}
 \end{equation}
 and
 \begin{equation}
    G(p, d)=\int_{0}^{\infty} y W_N(y, d) \exp(-py) dy.
\label{kwe71}
 \end{equation}

We suppose that $\bar k >0$ thus
$Re \ p > k$.
Therefore, we  can  investigate (\ref{ME}) in the domain
\begin{equation}
{\bf C}_{1/2, k}=\{p \in {\bf C}: \ p=\sqrt{k^2 +\lambda}+k,  \  Re \ \lambda > -1/2,  \  Re \ p > k \}.
\label{domp}
\end{equation}
Note that
\begin{equation} \label{lamviap}
Re \ p= k + \sqrt{k^2 + Re \lambda  + (Im \ p)^2}, 
\end{equation}
 this shows that in ${\bf C}_{1/2, k}$ we have $Re \ p > 2k -1/2$.
We can choose a polynomial $W_N$ such that 
\begin{equation}  \label{Gpd}
G(p,d)=p^{-2} (-1)^{N+1} \prod_{j=1}^N (\frac{1}{p} -  \frac{1}{2j + d_j}).
\end{equation}

Let us formulate an  auxiliary assertion.

\begin{lemma} \label{6.3} {
One has
\begin{equation}
 Re \ g_{\kappa}(p) \le - \min \{4 \exp(-4)  \kappa Re \ p, \quad 1/2   \}.
\label{g0}
\end{equation}  
}
\end{lemma}

{\bf Proof}. Estimate (\ref{g0})  follows from   (\ref{z0}) and (\ref{kweg}).
Indeed,  due to (\ref{Inteps}) we have
$$
Re \ g_{\kappa}(p) = \ \int_0^h \delta_{\kappa}(y-z_0)  (Re \ \exp(-py)-1)dy\le \int_0^h \delta_{\kappa} (\exp(-Re \ p y)- 1) dy
$$
that according to (\ref{z0})  gives $Re \ g_{\kappa}(p)  <  \exp( - 4 \kappa Re \ p )-1:=J(Re \ p) $.
Moreover,  
$
J( Re \ p)|  \le 4 \exp(-4) Re \ p| \kappa$  for $0 < Re \ p < \kappa^{-1}$ and  $J(Re \ p)  < 1/2$ for $|Re \ p | \ge \kappa^{-1}$.
\vspace{0.2cm}

Let us show that  in the case $ |p| > \kappa^{-3/4}$  equation (\ref{ME}) has no solutions with $Re \lambda >- c_0  \kappa$.

\begin{lemma} \label{plarge}
{If $|p| > \kappa^{-3/4}$ then for sufficiently small $\kappa$ solutions of (\ref{ME}) satisfy
\begin{equation}  \label{Repest}
Re \ p  < 2k  -c_1 k^2 \kappa^{1/4}
\end{equation}
For the corresponding $\lambda_k$ one has
\begin{equation} \label{lamest}
Re \ \lambda_k < -c_2 k^2 \kappa^{1/4}
\end{equation}
} 
\end{lemma}

{\bf Proof}. Relarion (\ref{lamviap}) shows  that  $Re \ p > Im \ p$,   thus $|p| > \kappa^{-3/4}$ entails  $4Re \ p >  \kappa^{-3/4}$.  
Then estimate (\ref{g0})  implies
$
 Re \ g_{\kappa}(p)    < -c_{2}\kappa^{1/4}.
$
Moreover,  
$
 \mu Re  \  G(p, d)    < -c_{3} \mu
$
and  for $Re \ p > -1/2$
\begin{equation} \label{tildeg}
 \tilde g_{\kappa}(p)  <  c_4 \mu |z_0| < c_5 \kappa^{5/3}.
\end{equation}
These  estimates  entail   that $Re \ S(p, k, d) < -c_4 \kappa^{1/4}$ and, 
therefore,
(\ref{Repest}) holds that, in turn,  by (\ref{lamviap}) gives us (\ref{lamest}).
\vspace{0.2cm}

Consider the case $|p| < \kappa^{-3/4}$.
Let us introduce a new unknown $\tilde p$ by
$
(2+\tilde p)k= p.
$
Then equation (\ref{ME}) can be  rewritten  as
\begin{equation}
\tilde p  =H(\tilde p, k),
\label{ME1}
 \end{equation}
where
\begin{equation}
H(\tilde p, k)= (\bar g_{\kappa}((2+\tilde p)k)  + \mu G( (2 +\tilde p)k, d))
\quad \label{Hpk}
 \end{equation}
and $\bar g_{\kappa}=g_{\kappa} + g_{\kappa}$.
Let us prove an  estimate of solutions to (\ref{ME1}).

 \begin{lemma} \label{6.4} { In the domain ${\mathcal D}_{\kappa, k}=\{ p: p \in {\bf C}_{1/2, k}, \ |p| < c_1\kappa^{-3/4} \}$ solutions of (\ref{ME1})
 satisfy
 \begin{equation}
 |\tilde p | < C_{1}\kappa^{1/4}.
 \label{tp}
 \end{equation}
 }
 \end{lemma}

{\bf Proof}.
 To prove this  lemma, we note that if $p \in {\mathcal D}_{\kappa, k}$ then 
\begin{equation}
|\bar g_{\kappa} (p)| < C_2\kappa^{1/4}, \quad |G(p, d)| < C_3.
\label{g0G}
\end{equation}
Moreover,  estimate (\ref{tildeg}) holds.
 Therefore, one obtains
\begin{equation} \label{Hunif}
|H(\tilde p, k)| < C_5 (\mu + \kappa^{1/4}).
\end{equation}
Now (\ref{Hpk}) and (\ref{ME1}) show that $\tilde p$  satisfies (\ref{tp}). The lemma is proved.

Let us consider  equation (\ref{ME1}).  Using the last lemma we note that, to resolve this equation,  
we can apply  a simple perturbation theory.
Relations (\ref{kwegt}), (\ref{kweg}),   (\ref{kwe70}) and (\ref{Gpd}) show that  in the  domain ${\mathcal D}_{\kappa, k}$ 
one has
$$
|\frac{\partial H(\tilde p, j)}{\partial \tilde p}|  < C_6  \mu.
$$
Now Lemma \ref{6.4} and
 the implicit function theorem  entail  that for sufficiently small $\kappa$ all roots $\tilde p$ of  eq.  (\ref{ME1}) lie in ${\mathcal D}_{\kappa, k}$ and can be found by 
contracting mappings.  For each fixed $k$  the solution $\tilde p_k$ of eq.  (\ref{ME1})  is unique in ${\mathcal D}_{\kappa, k}$.

\begin{lemma} \label{choiced}
{For sufficiently small $\kappa$ we can choose  $d_j \in (-1/2, 1/2)$, $j=1,..., N$, such that  for each $k \in \{1, ..., N \}$ eq. (\ref{ME1}) has a unique solution $\tilde p_k=0$ and for
$k > N$ any solution of eq. (\ref{ME1}) satisfy 
\begin{equation} \label{RekN}
Re \  \tilde p_k < - c\kappa^s, \quad s >0.
\end{equation}
} 
\end{lemma}

{\bf Proof}.   Due to Lemma \ref{plarge}, we can assume  $|p| < \kappa^{-3/4}$ and then, according to Lemma   \ref{6.4} for all $k$ solutions $\tilde p_k$ lie in the domain  defined by
inequality (\ref{tp}). 

Assume that $k \in \{1, ..., N \}$.
Let us fix $k$ and consider $p \in {\mathcal D}_{\kappa, k}$, i.e., $p$ close to $2k$. Then  for $\tilde p$ satisfying (\ref{tp}) from (\ref{Gpd}) we obtain the following asymptotic for $G(p, d)$: 
$$
G(p,d)= \mu (-1)^{N+1} (\bar a_k + \tilde a_k(\tilde p_k ,d)) (d_k - k \tilde p_k),
$$
where
$$
\bar a_k = \prod_{j=1, j \ne k}^{N+1} (\frac{1}{2k} -  \frac{1}{2j}),  
$$
and $\tilde a_k(\tilde p_k, d)$ is an analytic function such that $|\tilde a_k| =O(|\tilde p_k| + |d|)$ for small $\tilde p_k, d$.  
Eq. (\ref{ME1} takes the form
\begin{equation}\label{ME1d}
\tilde p_k= R_k(\tilde p_k, d, \kappa)  +    \mu (-1)^{M+1} \bar a_k d_k,
\end{equation}
where $R_k$ is an analytic function of $\tilde p_k$  and $d$ for   small $|\tilde p_k|, |d|$  such that 
$$
 \sup_{\tilde p_k, d: |d| < \mu^{1/2}, \tilde p_k \in  {\mathcal D}_{\kappa,k})} (|R_k| + |grad(R_k)|) < c \kappa^s  + c_2|b- b_c|,    
$$
for some $s >0$. 
Therefore, for small $b -b_c$  we can apply the Implicit Function Theorem to find need $d_k$ such that the root $\tilde p_k$ of (\ref{ME1d}) equal zero.

Let us consider the case $k >N$.
  We observe
that for $k > N$ and  $\tilde p$ satisfying (\ref{tp}) we have $\mu Re \ G < -c_2 \kappa^{2/3}$,  $Re \ \bar g_{\kappa}(p) < -c \kappa$,      
and thus $Re \ \tilde p <  - c_2 k \kappa^{2/3}$.  
The lemma is proved.
\vspace{0.2cm}

Finally,  we have obtained need estimates of solutions (\ref{kwc2}) for the case {\bf I}, $\beta_k=1$  and  $\nu=+ \infty$.
To finish our investigation of equation  (\ref{kwc2})  for the case of large $\nu$,  we compare equations (\ref{kwc2})  and its formal limit (\ref{kwc3}).  
We  assume that the function $U=V(y, d)$ is defined as above, by  (\ref{Uy}).

We observe that  for $\beta_k=1$ eq. (\ref{kwc2}) can be rewritten as
\begin{equation}
k \bar k^{-1}(\lambda)  \int_0^{+\infty}   y V(y,d) \exp(-(k +\bar k(\lambda)) y)dy=2 - R_k(\lambda, \nu, d),
\label{kwc3a}
\end{equation}
where 
$R_k=I_k + J_k$
and
\begin{equation}
I_k=k^2  \int_h^{+\infty}  \bar \psi_k(y,\lambda) \bar \rho_{\bar k} V_y(y,d) dy,
\label{kwI}
\end{equation}
\begin{equation}
J_k=-k^2  \int_0^{h}  ( \psi_k(y,\lambda) \rho_{\bar k} - \bar \psi_k(y) \bar \rho_{\bar k}) V_y(y, d) dy.
\label{kwJ}
\end{equation}
Under assumption (\ref{hviscos}), $\lambda << \nu^s$ for $s \in (0,1)$ and for sufficiently large $\nu$ the term $I_k$  satisfies the estimate
\begin{equation}
|I_k| < c_1 k^2 h^{N+3} \exp(-kh) < c_2 \nu^{-4}, 
\label{kwest}
\end{equation}
which is uniform in $k$.
To estimate $J_k$ we use  the inequality 
\begin{equation}
|J_k| \le k^2  \int_0^{h}  (| (\psi_k(y,\lambda) -\bar \psi_k(y)) \rho_{\bar k}| + |\bar \psi_k(y) (\bar \rho_{\bar k}(y) -\rho_k(y))|) |V_y(y,d)| dy.
\label{kwJ2}
\end{equation}
Now we apply (\ref{psiFa}), definition  (\ref{Uy}) of $V_y$, estimates 
 (\ref{psiFa1}), (\ref{tildepsi}) and (\ref{rhoest1}) that gives
$$
\sup_{y \in [z_0, h]}  |V_y(y, d)| <  c_3 (\kappa^{-2}  + h^{N+3})\exp(- k \kappa).  
$$
Then we see that  
$$
|J_k|  <  c_4 k^2 (\kappa^{-2}  + h^{N+3}) (k^2 \bar k^{-1} \exp(-\bar kh) + k \bar k^{-1} \lambda \nu^{-1})  \exp(- k \kappa) <
$$
$$
 < c_5 (1 + \kappa^{-2}) \nu^{-1/2}  
$$
for some  $c_5 >0$ and sufficiently large $\nu > \nu_0(\kappa)$.
 Note that in this estimate the constant $c_5$ is uniform in $k$.  
We obtain finally the uniform in $k$ estimate
\begin{equation}
|R_k(\lambda, \nu, d)|   <  c_6  \nu^{-1/2}.
\label{RKK}
\end{equation}

Therefore, the  analysis of equation (\ref{kwc3a}) can be made by the same arguments as above in the case of formal limit $\nu \to +\infty$ 
that allows us to prove the lemma.

\begin{lemma} \label{LL}
{ Let  assumptions (\ref{hviscos}) hold,
 $N$ be a positive integer  and  $V(y, d)$ is defined by  
(\ref{Uy}).   Moreover, let $|\lambda| < \nu^{3/4}$.  Then  we can choose such parameters $\kappa$ and $d$  
in  relation (\ref{Uy}  that for sufficiently large $\nu$ the roots  
 $\lambda(k, \nu)$ of equation (\ref{kwc3a}), which lie   in the domain $ |\lambda| < \nu^{3/4}$,  satisfy
\begin{equation}
\lambda(k, \nu)=0     \quad  k =1, ..., N
\label{Spec0A}
\end{equation}
and  
\begin{equation}
  Re \ \lambda(k,  \nu) < - \delta_N  \quad  k > N,
  \label{SpecA}
\end{equation}
where a positive $\delta_N$ is uniform in $k$ and in $\nu$ as  $\nu \to \infty$.
}
\end{lemma}

{\bf Proof}. 
We repeat the proof of Lemma \ref{choiced}.  Since  $R_k(\lambda, \nu, d)$ satisfies  estimate (\ref{choiced}),   we obtain new $d_j=\bar d_j(\nu)$, which are small  
perturbations of $d_j$ obtained  in Lemma \ref{choiced}.  We have  
$\bar d_j(\nu) -d_j=\tilde d_j(\nu)$, where   $\tilde d_j(\nu)\to 0$ as $\nu \to +\infty$. Therefore, the $\sup |U_y(y, \nu)-  V_y(y,d)| \to 0$ as $\nu \to +\infty$ 
and (\ref{Spec0A}) holds.  Inequalities  (\ref{Spec0A}) are fulfilled for sufficiently large $\nu$.  It follows from equation (\ref{kwc3a}) and  estimate (\ref{choiced}), which is uniform in $k$.   
The lemma is proved.

\vspace{0.2cm}

{\em  Case {\bf II}}.

Let us consider now the second case {\bf II}.
Relations 
(\ref{psiF2}) and (\ref{rhoest1}) imply that 
$$
|\psi_k(y, \lambda)| <  c_1 \nu/|\lambda|,   \quad |\rho_{\bar k}(y) |  < c_2 \exp(-\bar k  y) |\bar k|^{-1}, 
$$
where $|\bar k| >  |\lambda|^{1/2}$.  Moreover, $|y U_y| < C h^{M+2}$. 
Thus the left hand side of eq. (\ref{kwc2}) is not more than  $R=C\nu  \log(\nu)^{M+2} k \bar |k|^{-1} |\lambda|^{-3/2}$. 
For $\nu \to +\infty$ and $|\lambda > \nu^{3/4}$ one has $R < c\nu^{-1/8}$.  Therefore,  eq. (\ref{kwc2}) has no solutions in the case {\bf II}.

The assertion of Proposition \ref{6.2} follows from this lemma and Lemma \ref{LL}.

\subsection{ Eigenfunctions of $L$ with zero eigenvalues}
\label{sec8.4}

Let us consider the eigenfunctions $e_k$ of $L$ with the zero eigenvalues. We have $2N+1$  eigenfunctions  including  the trivial one 
$e_0=(0, 1)$.  All the rest eigenfunctions   have the form
\begin{equation}
e_k=\exp(ikx)(\omega_k, \theta_k)^{tr},   \quad e_{-k}=\exp(-ikx)(\omega_{-k}, \theta_{-k})^{tr}, 
\label{eig1}
\end{equation}
where $k=1,2, ..., N$ and 
\begin{equation}
\omega_k= i k A_k  \frac{\sinh (k(h-y))}{\sinh(k h)},  \quad i=\sqrt{-1},
\label{eig2}
\end{equation}
\begin{equation}
\theta_k= A_k \Theta_{k}(y),
\label{eig3}
\end{equation}
where $A_k(\nu)$ are  constants. 
The functions $\Theta_k$ are defined by
\begin{equation}
\frac{d^2\Theta_{k}(y)}{dy^2} -  k^2 \Theta_k= i k \Psi_k(y)U_y(y),
\label{eig4}
\end{equation}
where $\Psi_k$ is given by (\ref{psiF4}) with $\lambda=0$ and $w_k(0,\lambda)=1$.
Thus,  $\Psi_k$  is defined 
\begin{equation}
\Psi_k(y)= i\frac{ y\sinh(kh) \cosh(k(h-y)) - h \sinh(ky) }{2\sinh^2( k h)}.    
\label{asp}
\end{equation}
Asymptotics of $\omega_k$ is given by (\ref{psiF}).
Using relations for $U$, estimates from the previous subsection  and the definition of $\mu, \kappa, z_0, z_1$,
one has the following asymptotics
\begin{equation}
\Theta_k(y) = k^{-1} \exp(-k y) (1 + \tilde r_k(y, \nu)), \quad y \in (z_0 ,h),
\label{asttheta2}
\end{equation}
where
\begin{equation}
\sup_{y \in [z_0, h]}|\tilde r_k(y, \nu)| < C_1 \nu^{-1/100}.
\label{asttheta2r}
\end{equation}

In coming subsection  we  construct the operator $L^*$, formally conjugate to $L$, and
the corresponding eigenfunctions.

\subsection{Conjugate spectral problem }
\label{sec8.5}

Consider spectral problem for conjugate operator $L^*$.  We follow \cite{Sengul} but use the vortex-stream formulation.  

Let us denote by $\rho$ and $\tilde \rho$ the pairs $\rho=(\omega, u)^{tr}$
and $\tilde \rho=(z, v)^{tr}$. Let us calculate the quadratic form $( L\rho, \tilde \rho ) $, where
$(\rho, \tilde \rho)$ denotes a natural inner product defined  by
$$
(\rho, \tilde \rho)=\langle \omega, \tilde \omega \rangle + \langle u, \tilde u \rangle.
$$
We obtain  the  relations
\begin{equation}
\langle \Delta u + \Psi_x U_y, v\rangle =
 \langle u, \Delta v \rangle - \langle \omega, \Phi(z)_x\rangle + J,
\label{conja}
\end{equation}
where
$$
J= \int_0^{2\pi} (u_y(x,y) v(x,y) - u(x,y) v_y(x,y)) dx)\vert_{y=0}^{y=h},
$$
and  $\Phi$ is  a unique solution
of the boundary value problem
\begin{equation}
\Delta \Phi = U_y v,
\label{conj1}
\end{equation}
\begin{equation}
\Phi(x,0)=\Phi(x, h)=0.
\label{conj2}
\end{equation}
Note that
\begin{equation}
\langle \nu\Delta \omega, z\rangle =
 \langle \omega, \nu \Delta z\rangle + \nu \int_0^{2\pi} (\omega_y(x,y) z(x,y) -
 z_y(x,y) \omega(x,y)  dx)\vert_{y=0}^{y=h}.
\label{conjb}
\end{equation}
Let us compute the eigenfunctions of the conjugate operator $L^*$.
Note that these eigenfunctions have been found in \cite{Sengul} in the three-dimensional case for linear profiles $U(y)$.

Using  boundary conditions
(\ref{boundstream2}) and (\ref{Maran1}) (\ref{Maran2}), by (\ref{conja}) and (\ref{conjb})
we obtain that the discrete spectrum of the operator $L^*$ is defined by
 the following equations
\begin{equation}
\lambda z = \nu \Delta z - \Phi(v)_x,
\label{conj4}
\end{equation}
\begin{equation}
\lambda v = \Delta v,
\label{conj3}
\end{equation}
 under the boundary conditions
\begin{equation}
v_y(x,y)\vert_{y=h}=0, \quad  v_y(x,y)\vert_{y=0}= \nu z_{xy}(x,y)\vert_{y=0},
\label{conj5}
\end{equation}
\begin{equation}
z(x,0)=z(x, h)=0.
\label{conj6}
\end{equation}

Now the eigenfunctions of $L^*$  with the zero eigenvalues can be found.  We obtain
\begin{equation}
\tilde e_k=\exp(ikx)(z_k, \tilde\theta_k)^{tr}, \quad k \in \{-N,  ..., N\} ,
\label{eig1m}
\end{equation}
where
\begin{equation}
\tilde\theta_k= \tilde A_k \bar \theta_k(y), \quad \bar\theta_k=\frac{\cosh( k(h-y))}{\sinh( k h)},
\label{eig3m}
\end{equation}
\begin{equation}
z_k=  \nu^{-1}  \zeta_{k}(y).
\label{eig2m}
\end{equation}
Here
 $\zeta_k$ are defined as  solutions of the boundary value problem
 $$
 \frac{d^2\zeta_k}{dy^2} -k^2 \zeta_k= - ikU(y)\Phi_k(y)
 $$
 $$
 \zeta_k(h)=0,  \quad \frac{d\zeta_k(y)}{dy}\vert_{y=0}=i\tilde A_k,
 $$
 where  $\Phi_k(y)$ are defined by
 $$
 \frac{d^2\Phi_k}{dy^2} -k^2 \Phi_k=  \tilde A_k U(y) \bar \theta_k(y), \quad \Phi_k(0)=\Phi_k(h)=0.
 $$
  We have the  estimate  
 \begin{equation} \label{estz}
 |z_k|_{C^2(\Omega)} < C_0\nu^{-1}.
\end{equation}
The exact form of the functions $z_k$ is not essential.
Only expressions for
the functions $\psi_k$ and $\theta_k, \tilde \theta_k$ and estimate (\ref{estz}) are involved in the further statement.

Relations (\ref{eig3m}), (\ref{asttheta2}) and (\ref{estz}) show that one can  adjust
the constants $A_k$ and $\tilde A_k$  such that 
\begin{equation}
\langle e_l, \tilde e_j \rangle=\delta_{lj}.
\label{eig6m}
\end{equation}
Note that $\tilde e_0=\tilde A_0 (0,1)^{tr}$.  

To obtain real value eigenfunctions, we  take real and imaginary parts of these complex eigenfunctions.
The real parts of the eigenfunctions, where $ \omega_k, \theta_k$ are proportional to $\sin(k x), \ \cos(k x)$ respectively,
 are enabled by the upper index $+$,  and the imaginary parts, where $ \omega_k, \theta_k$ are proportional to $\cos(k x), \ \sin(k x)$,
are denoted by the upper index $-$.  The  real eigenfunctions of $L$ have the form

\begin{equation}
e_k^+=(\omega_k(y) \sin(kx), \theta_k(y) \cos(kx))^{tr}, 
\end{equation}
\begin{equation}
 e_k^-=(-\omega_k(y) \cos(kx) , \theta_k(y) \sin(kx))^{tr}.
\label{eigpm}
\end{equation}
 Respectively, the  real eigenfunctions of $L^*$ are
\begin{equation}
\tilde e_k^+=(\tilde \omega_k(y) \sin(kx), \tilde \theta_k(y) \cos(kx))^{tr}, 
\end{equation}
\begin{equation}
 \tilde  e_k^-=(-\tilde \omega_k(y) \cos(kx) ,\tilde  \theta_k(y) \sin(kx))^{tr}.
\label{eigpms}
\end{equation}
We  have the  relations
\begin{equation}
\tilde\theta_k^{+}=a_k  \cosh( k(h-y)), \quad
\label{eigplus}
\end{equation}
where $a_k$ are coefficients. 

The next lemma 
concludes the investigation of spectral properties of the operator $L$.

\begin{lemma} \label{simplespec}
{For sufficiently large $\nu$ the eigenvalue $0$ of the operator $L$ has  the multiplicity $N$. The eigenvalue $0$ has no generalized eigenfunctions.}
\end{lemma}

{\bf Proof}. Let us check that generalized eigenfunctions are absent.
Since $0$ has a finite multiplicity, we can use the Jordan representation.  Assume that there exists a generalized eigenfunction  $e_g$.
Then
$
L e_g=b=\sum_{j=1}^N b_l e_l
$
for some $b_l$, $ l \in \{1, ..., N\}$,  where $b \ne 0$.  Then  $\langle b, \tilde e_k \rangle=0$ for all $k \in \{1,..., N\}$. Eigenfunctions $\tilde e_k$ and $e_l$ are biorthogonal according to  (\ref{eig6m}).
This implies that all coefficients $b_l=0$ and the lemma is proved.      
\mbox

\subsection{ Estimates for semigroup $\exp(Lt)$ }
\label{sec: 8.6}

The operator $L$ is sectorial and, according to  Lemma \ref{6.2},   satisfies  the Spectral Gap Condition. Therefore \cite{He}
for some $\beta(\nu) > 0$ and for all $v=(\tilde \omega, w)^{tr}$ such that
$$
\langle v, \tilde e_j^{\pm} \rangle =0, \quad j=1, ..., N,
$$
and
$$
\langle v, \tilde e_0 \rangle=0,  \quad \tilde e_0=(0, 1).
$$
one has the following estimates 
\begin{equation}
|| \exp(Lt) v|| \le C_0(\nu) \exp(-  \beta(\nu)t) ||v||,  \quad
\label {semigroup1}
\end{equation}
\begin{equation}
|| \exp(Lt) v||_{{\mathcal H}_\alpha} \le C_1(\nu) b_{\alpha}(t) \exp(- \beta(\nu)t) ||v||
\label {semigroup2}
\end{equation}
Here
$$
b_{\alpha}(t) = t^{-\alpha}, \quad 0 < t \le \beta^{-1}(\nu)
$$
and
$$
b_{\alpha}(t) = \beta^{-\alpha}, \quad  t > \beta^{-1}(\nu).
$$
Estimates (\ref{semigroup1}) and (\ref{semigroup2}) are important in the proof of  existence of the invariant manifold.

\section{Finite dimensional invariant manifold
}
\label{sec: 9}

 Assume
$\gamma >0$ is a small parameter. 
 In this section, we reduce the Navier -Stokes  dynamics to a system of ordinary differential equations 
  following  Section 4.  Let  $E_N$ be
the finite dimensional subspace
$E_N=Span \{e_0, e_1^+,  ..., e_{N}^+, e_1^-, ..., e_N^{-} \}$
of the phase space ${\mathcal H}$, where $e_j^{\pm}=(\omega_j^{\pm}, \theta_j^{\pm})^{tr}$ are the eigenfunctions of the operator
$L$ with the  zero eigenvalues.
Let ${\bf P}_N$ be a  projection operator
on $E_N$ and ${\bf Q}_N={\bf I} - {\bf P}_N$.
The components of ${\bf P}_N$ are defined by
\begin{equation}
{\bf P}_{1, N} v= \sum_{j=1}^{N} \langle \tilde \omega, \tilde \omega_j^{+} \rangle \omega_j^{+} +   \sum_{j=1}^{N}\langle \tilde \omega, \tilde \omega_j^{-} \rangle \omega_j^{-},
\label{Pr2}
\end{equation}
\begin{equation}
{\bf P}_{2, N} v =(2\pi h)^{-1} \langle w, 1\rangle +  \sum_{j=1}^{N} \langle w, \tilde \theta_j ^+ \rangle \theta_j^+   +\sum_{j=1}^{N} \langle w, \tilde \theta_j ^- \rangle \theta_j^-,
\label{Pr2a}
\end{equation}
where  $v=(\tilde \omega, w)^{tr}$ and  $\tilde e_j^{\pm} =(\tilde \omega_i^{\pm}, \tilde \theta_i^{\pm})$ are eigenfunctions of the conjugate operator $L^*$ with zero eigenvalues $\lambda=0$
found in Sect. \ref{sec8.5}.
Let us rewrite system (\ref{eveq10}), (\ref{eveq11}) as
\begin{equation}
  \omega_t=\nu \Delta  \omega   - \{\psi, \omega \},
\label{eveq10a}
\end{equation}
\begin{equation}
w_t= \Delta w   - \{\psi, U +\gamma u_1 + w \} + \gamma^2 \eta_1.
\label{eveq11b}
\end{equation}

First we transform equations (\ref{eveq10a})-(\ref{eveq11b}) to a standard system with "fast" and "slow" modes.
We follows  Section 4. Let us introduce auxiliary functions $R_{\omega}(X)$, $R_{\psi}(X)$ and $R_{w}(X)$ by
$$
 R_{\omega}(X) =\sum_{j=1}^{N} X_j^+ \omega_j^+ + \sum_{j=1}^{N} X_j^- \omega_j^-,  \quad  R_{\psi}(X)=\sum_{j=0}^{N} X_j^ + \psi_j^+ +   \sum_{j=1}^{N} X_j^ -  \psi_j^-, 
$$
$$
 R_w(X)=X_0 +  \sum_{j=1}^{N} X_j^+ \theta_j^+    +  \sum_{j=1}^{N} X_j^- \theta_j^- ,
$$
and represent $ \omega, \psi$ and $w$ by
\begin{equation}
\omega= \gamma R_{\omega}(X) + \hat \omega,  \quad 
 \quad \psi= \gamma R_{\psi}(X) +\hat \psi,   
\label{om2}
\end{equation}
\begin{equation}
w= \gamma R_w(X) +\hat w,  
\label{w2}
\end{equation}
where
$
{\bf P}_N (\hat \omega, \hat w)^{tr}=0,
$
$X_i^{\pm}(t)$ are unknown functions,
 $X=(X_0, X_1^+, ..., X_{N}^+, X_1^-,..., X_N^{-})^{tr}$.  

We  substitute relations (\ref{om2})-(\ref{w2}) in eqs. (\ref{eveq10a}) and (\ref{eveq11b}).
 After some transformations (following Section 4)  one obtains  the system
\begin{equation}
\frac{dX_i^{\pm}}{dt}= \gamma  (G_i^{\pm}( X) +  M_{i}^{\pm}( X) + f_i^{\pm} + F_i^{\pm}( X, \hat \omega,\hat w, \gamma)),
\label{X2}
\end{equation}
\begin{equation}
\hat \omega_t=\nu \Delta \hat \omega +  {\bf P}_{1,N} F(X, \hat \omega,\hat w, \gamma),
\label{hom2}
\end{equation}
\begin{equation}
\hat w_t= \Delta \hat w   - \{\hat \psi, U\} + {\bf P}_{2,N} G( X, \hat \omega,\hat w, \gamma),
\label{hw2}
\end{equation}
where
in eqs.(\ref{hom2}) and (\ref{hw2}) 
\begin{equation}
F=\{\gamma R_{\psi}(X) + \hat \psi, \gamma R_{\omega}(X) + \hat \omega\},
\label{F2}
\end{equation}
\begin{equation}
G=\{\gamma R_{\psi}(X) + \hat \psi, \gamma R_{w}(X) + \gamma u_1 + \hat w\} 
+\gamma^2 \eta_1.
\label{G2}
\end{equation}
The functions $G_i^{\pm}(X)$ and $M_i^{\pm}(X)$ give main contributions in the right hand sides of eqs.   (\ref{X2}) and $F_i^{\pm}$ are small corrections for  small $\gamma$.  One has 
$$
G_i^{\pm}(X)=\langle\{R_{\psi}(X), R_{w}(X)\}, \tilde \theta_i^{\pm} \rangle +  \langle\{R_{\psi}(X), R_{\omega}(X)\}, \tilde \omega_i^{\pm} \rangle,    
$$
$$
M_i^{\pm}(X)=\langle\{R_{\psi}(X), u_1\}, \tilde \theta_i^{\pm} \rangle.
$$
These terms can be rewritten in a more explicit form as
\begin{equation}\label{Gpm}
  G_i^{+}( X) =\sum_{j,l=1}^N G_{ijl}^{+++} X_j^+ X_l^+  +  G_{ijl}^{--+} X_j^- X_l^-,  
\end{equation}
\begin{equation}\label{Gpm-}
    G_i^{-}( X) = \sum_{j,l=1}^N G_{ijl}^{+--} X_j^+ X_l^-
\end{equation}
and
\begin{equation}\label{Mpm}
  M_i^{+}( X) = \sum_{j=1}^N M_{ij}^{++} X_j^+
   + \sum_{j=1}^N M_{ij}^{+-} X_j^- , 
\end{equation}
\begin{equation}\label{Mpm-}
 M_i^{-}( X) =\sum_{j=1}^N M_{ij}^{-+} X_j^+
   + \sum_{j=1}^N M_{ij}^{--} X_j^-  .
\end{equation}
 Note that in eqs. (\ref{Gpm}) - (\ref{Mpm-}) 
 all the rest possible  terms vanish  since they are defined by integrals over $x$ of functions  odd in $x$

The coefficients in (\ref{Gpm}),(\ref{Gpm-}), (\ref{Mpm}) and (\ref{Mpm-}) are defined by
\begin{equation}
   M_{ij}^{\pm \pm}(u_1)=
\langle \{\psi_j^{\pm}, \tilde \theta_i^{\pm}  \},  u_1\rangle,
\label{MatrM}
\end{equation}
\begin{equation}
   G_{ijl}^{+ + +}=\langle \{\psi_j^{+}, \theta_l^{+}\}, \tilde \theta_i^{+} \rangle + O(\nu^{-1})
\label{G+}
\end{equation}
\begin{equation}
   G_{ijl}^{- - +}=\langle \{\psi_j^{-}, \theta_l^{-}\}, \tilde \theta_i^{+} \rangle + O(\nu^{-1})
\label{G-}
\end{equation}
\begin{equation}
   G_{ijl}^{ + --}=\langle \{\psi_j^{+}, \theta_l^{-}\} + \{\psi_l^{-}, \theta_j^{+}\}, \tilde \theta_i^{-} \rangle + O(\nu^{-1})
\label{matrG}
\end{equation}
for large $\nu$.   Here we have used estimate (\ref{estz}), which implies that the terms,  where  $\tilde \omega_j$ are involved, have the order $O(\nu^{-1})$.  One has
\begin{equation}
  f_{i}^{\pm}=
\langle \eta_1,  \tilde \theta_i^{\pm} \rangle.
\label{fipm}
\end{equation}
The  terms  $F_i^{\pm}$ are  defined by 
\begin{equation}
F_i^{\pm}= \gamma^{-1} (F_i^{\pm, \omega} + F_i^{\pm, w}) , 
\label{Fi2}
\end{equation}
where
$$
F_i^{\pm, \omega}=\langle \{  \gamma R_{\psi}(X), \hat \omega \} + \{ \hat \psi, \gamma R_{\omega}(X)  + \hat \omega \}, \tilde \omega_i^{\pm} \rangle, 
$$
$$
F_i^{\pm, w}=\langle  \{  \gamma R_{\psi}(X), \hat w \}   + \{\hat \psi, \gamma R_{w}(X)  + \gamma u_1 + \hat w \}, \tilde \theta_i^{\pm} \rangle). 
$$
We consider   equations (\ref{X2}), (\ref{hom2}) and (\ref{hw2})  in the domain
\begin{equation}
  {\mathcal D}_{\gamma, R_0, C_1, C_2, \alpha} =\{(X, \hat w, \hat {\omega}): \ |X| <  R_0, \ ||\hat \omega||_{\alpha} < C_1\gamma^{3/2},
   ||\hat w||_{1 +\alpha}  < C_2\gamma^{3/2} \},
\label{Dom}
\end{equation}
where $\alpha > 3/4$.
 Let us define the vector field $V$  on the  ball ${\mathcal B}(R_0)^{2N} \subset {\bf R}^{2N}$ by
$$V(X)=(V_1^+(X), ..., V_N^+(X),  V_1^-(X), ..., V_N^-(X))$$.

\begin{lemma} \label{7.1}
{
Let $r  \in (0,1)$ and $\alpha \in (3/4, 1)$.
Assume $\gamma > 0$ is small enough:
$
\gamma < \gamma_0(N, \nu, R_0, r, \alpha).
$
Then the local semiflow $S^t$, defined by equations  (\ref{X2}),(\ref{hom2}), and (\ref{hw2}) has a locally invariant in the set ${\mathcal D}_{\gamma, R_0, C_1, C_2, \alpha} \subset \mathcal H$  
and locally attracting manifold ${\mathcal M}_{2N+1,\gamma}$. This manifold is  defined by
\begin{equation}
 \hat \omega= \hat \omega_0(X, \gamma),   \quad \hat w= 
   \hat w_0(X, \gamma),
\label{rW}
\end{equation}
where
$\hat \omega_0(X, \gamma)$, $\hat w_0(X,\gamma)$ are maps from the ball ${\mathcal B}^{2N+1}(R_0)=\{X: |X| <  R_0\}$ to $H_{\alpha}$ and $H_{1+\alpha}$ respectively,
bounded in $C^{1+r}$ -norm :
\begin{equation}
  |\hat \omega_0(X, \gamma)|_{C^{1+r}({\mathcal B}^{2N+1}(R_0))} < C_3\gamma^2, 
\label{rWe}
\end{equation}
\begin{equation}
   |\hat w_0(X, \gamma)|_{C^{1+r}({\mathcal B}^{2N+1}(R_0))} < C_4 \gamma^2.
\label{rWe1}
\end{equation}
The restriction of  the semiflow $S^t$ on ${\mathcal M}_{2N+1, \gamma}$  is defined by the vector field $V(X)$. The corresponding differential equations
have the form
\begin{equation}
  \frac{dX_i^{\pm}}{dt}=\gamma (V_{i}^{\pm}(X)  +  \phi_i^{\pm}(X, \gamma)), 
\label{maineq1}
\end{equation}
where  $X=(X^{+}, X^{-})$,  
\begin{equation}  \label{Vpm}
V_i^{\pm}(X)=G_i^{\pm}( X) +  M_{i}^{\pm}( X) + f_i^{\pm}) 
\end{equation}
\begin{equation}
  \frac{dX_0}{dt}=0,
\label{maineq0}
\end{equation}
and the corrections $\phi_i^{\pm}(X, \gamma)=
 F_i^{\pm}( X, \hat \omega_0(X,\gamma),\hat w_0(X, \gamma), \gamma)
$ satisfy the estimates
\begin{equation}
  |\phi_i^{\pm}|, |D_X \phi_i^{\pm}|  < c_1 \gamma^s,   \quad  s>0.
\label{maineq4}
\end{equation}
}
\end{lemma}

This assertion is a consequence of Lemma \ref{rvf}.
In coming sections we investigate system  (\ref{maineq1}).

\section{Quadratic systems}
\label{sec:5}

For sufficiently small $\gamma$  we can remove  small corrections $\phi_i^{\pm}$ in the right hands of  (\ref{maineq1}).    Then we obtain a  system of differential equations with quadratic nonlinearities.   
Let us consider a general class of such quadratic systems
\begin{equation}
   \frac{dX}{dt}=   { K}(X)    + { M} X + g,
\label{gks}
\end{equation}
where  $X=(X_1, ..., X_N), \ { K}=(K_1, ...,K_N),  { g}=(g_1, ..., g_N) \in {\bf R}^N$, $K(X)$ is a quadratic term
defined by
$$
{K}_i(X)= \sum_{j=1}^N \sum_{l=1}^N K_{ijl} X_j X_l,
$$
and $MX$ is a linear operator
$$
({ M}X)_i= \sum_{j=1}^N M_{ij} X_j.
$$
System (\ref{gks}) defines a local semiflow $S^t(g,{ M})$ in the ball ${\mathcal B}^N(R_0) \subset {\bf R}^N$ of the radius $R_0$
centered at $0$. We shall consider the vector $g$ and the matrix $ M$ as parameters of this semiflow whereas
the entries $K_{ijl}$ will be fixed.

Let us formulate an assumption on entries $K_{ijl}$.
We  present  $X$ as a pair $X=(Y, Z)$, where
$$
Y_l=X_{l},  \quad l \in I_p, \quad Z_j=X_{j+p}, \quad l \in J_p.
$$
Here  $I_p=\{1,..., p \}$ and $J_p=\{p+1, ..., N \}$ are  subsets of $\{1,..., N\}$.
Then the system (\ref{gks}) can be rewritten as
 \begin{equation}
   \frac{dY}{dt}=   { K}^{(1)} (Y) + { K}^{(2)} (Y, Z) + { K}^{(3)}(Z) + { R} Y + { P} Z
 +  f,
\label{gks2}
\end{equation}
 \begin{equation}
   \frac{dZ}{dt}=  \tilde { K}^{(1)} (Y) + \tilde{ K}^{(2)} (Y, Z) + \tilde { K}^{(3)} (Z)  +  \tilde { R} Y +  \tilde { P} Z
 +  \tilde f,
\label{gks3}
\end{equation}
where  for $i=1,..., p$ 
\begin{equation}
   { K}^{(1)}_i(Y)=  \sum_{j \in I_p} \sum_{l \in I_p}  K_{ijl}^{(1)} Y_{j} Y_{l}, \quad { K}^{(3)}_i(Z)=  \sum_{j \in J_p} \sum_{l \in J_p}  K_{ijl}^{(3)} Z_{j} Z_{l},
\label{KY1}
\end{equation}
\begin{equation}
   { K}^{(2)}_i(Y,Z)=  \sum_{j \in I_p} \sum_{l \in J_p}  K_{ijl}^{(2)} Y_{j} Z_{l}, 
\label{KY2}
\end{equation}
and for $k=1,..., N-p$
\begin{equation}
    \tilde { K}_k^{(1)}(Y)= \sum_{j\in I_p} \sum_{l \in I_p}  \tilde K_{kjl}^{(1)} Y_{j} Y_{l}, \quad \tilde { K}_k^{(3)}(Z)=\sum_{j \in J_p} \sum_{l \in J_p}  \tilde  K_{kjl}^{(3)} Z_{j} Z_{l},
\label{KY3}
\end{equation}
\begin{equation}  \label{K2YZ}
\tilde { K}_k^{(2)}(Y,Z)= \sum_{j \in I_p} \sum_{l \in J_p}  \tilde  K_{kjl}^{(2)} Y_{j} Z_{l}.
\end{equation}
The linear terms $MX$ take the form 
 \begin{equation}
 ({ R} Y)_i= \sum_{j \in I_p} R_{ij} Y_j,  \quad  (\tilde { R} Y)_k= \sum_{j \in I_p} \tilde R_{kj} Y_j,
\label{gks3c}
\end{equation}
\begin{equation}
  ( { P} Z)_i=\sum_{j \in J_p}  P_{ij} {Z_j},
 \quad  ({\tilde  P} Z)_k= \sum_{j \in J_p} \tilde P_{kj} Z_j,
\label{gks3d}
\end{equation}
and $f=(f_{1}, ..., f_{p}), \ \tilde f=(\tilde f_1, ..., \tilde f_{N-p})$.
We denote by  $S^t({\mathcal P})$ the local semiflow defined by
(\ref{gks2}) and (\ref{gks3}). Here ${\mathcal P}$ is a semiflow parameter, $ {\mathcal P}=\{f, \tilde f, {P}, \tilde { P},
{R}, \tilde { R} \}$.  Let us formulate an assumption on quadratic terms $K_i(X)$. 
\vspace{0.2cm}

{\bf $p$-Decomposition Condition} \label{ass40}
{\em
Suppose  entries $K_{ijl}$ satisfy  the following condition.
For some $p$ there exists a decomposition $Z=(X, Y)$, where $X \in {\bf R}^p$ and $Y \in {\bf R}^{N-p}$ such that
 the  linear system  
\begin{equation}
	     \sum_{i \in J_p} \tilde K_{ijl}^{(1)}  u_i= b_{jl}, \quad l, j \in I_p
\label{barK}
\end{equation}
has a solution $u$   for all   $b_{jl}$.
}
\vspace{0.2cm}

Clearly that for $N > p^2 + p$ and  generic matrices $K$ this condition is valid.

Let us formulate  some conditions to  the matrices ${ R}, \tilde { R}, {P}$ and $\tilde { P}$.
Let $\epsilon >0$ be  a  parameter.  We suppose that
\begin{equation}
  \tilde P_{ij}=-\xi^{-1} \delta_{ij},    \quad  i=1,..., N-p, \ j=1, ...,
\label{gks4c}
\end{equation}
where $\delta_{ij}$ is the Kronecker symbol,
\begin{equation}
  \tilde R_{ij}=0, \quad  \tilde f_i=0, \quad i=1,..., N-p, \ j=1,...,p,
\label{gks4d}
\end{equation}
\begin{equation}
   P_{ij}=\xi^{-1} T_{ij}, \quad |T_{ij}| < C_0, \quad i=1,..., p, \ j=1,..., N-p,
\label{gks5c}
\end{equation}
\begin{equation}
 |R_{ij}| < C, \quad i=1,..., p, \ j =1,..., p,
\label{gks5d}
\end{equation}

\begin{lemma} \label{quadr1} {Assume  (\ref{gks4c}),  (\ref{gks4d}), (\ref{gks5c}) and (\ref{gks5d})  hold and $R_0 >0$.
For sufficiently small positive $\xi < \xi_0(R_0, r, f)$ the local semiflow $S^t({\mathcal P})$ defined by  system (\ref{gks2}), (\ref{gks3}) 
has a locally invariant in the domain
\begin{equation} \label{domainR0}
{\mathcal  D}_{R_0}=\{X:   \quad  |Y| < R_0 \} 
\end{equation}
 and locally attracting manifold ${\mathcal M}_{\mathcal P}$.  This manifold is  defined by equations
\begin{equation}
  Z= \xi (\tilde { K}^{(1)}(Y)  +  W(Y,\xi)), \quad Y \in {\mathcal B}^p(R_0)
\label{gks6c}
\end{equation}
where $W$ is a $C^1$ smooth map defined on the ball ${\mathcal B}^p(R_0)$ to ${\bf R}^{N-p}$ and such that
\begin{equation}
|W(\cdot, \xi)|_{C^1({\mathcal B}^p(R_0))}  < C_1\xi^s, \quad s > 0.
\label{gks7c}
\end{equation}}
\end{lemma}

{\bf Proof}. Let us introduce a new variable $w$ by
\begin{equation}
  Z= \xi (\tilde { K}^{(1)}(Y)  +  w).
\label{gks8c}
\end{equation}
and the rescaled time by $t=\xi\tau$.
Then for $Y, w$ one obtains the following system
\begin{equation}
   \frac{dY}{d\tau}=   \xi G(Y, w, \xi),
\label{gk1}
\end{equation}
 \begin{equation}
   \frac{dw}{d\tau}=  \xi F(Y, w, \xi)
    -  w,
\label{gk2}
\end{equation}
where
$$
    G(Y, w, \xi)={ K}^{(1)} (Y) + \xi {K}^{(2)} (Y, \tilde { K}^{(1)}(Y)  +  w) +
$$
$$
    +
  \xi^2 { K}^{(3)}(\tilde { K}^{(1)}(Y)  +  w) + { R} Y + { T}(\tilde { K}^{(1)}(Y)  +  w)
 +  f,
$$
$$
    F(Y, w, \xi)=   \tilde{ K}^{(2)} (Y, \tilde { K}^{(1)}(Y)  +  w)+ $$
 $$   +
 \xi  \tilde { K}^{(3)} (\tilde { K}^{(1)}(Y)  +  w) +  h(Y, w, \xi),
$$
\begin{equation}
    h(Y, w, \xi)= - (D_Y \tilde { K}^{(1)}(Y))  G(Y, w, \xi).
\label{gk4}
\end{equation}

 Equations (\ref{gk1}), (\ref{gk2}) form  a typical system involving slow ($Y$) and fast ($w$) variables. 
Existence of a locally invariant manifold for  this system can be shown by the well known results 
(see \cite{He,CLu, CFNT, Bates,Tem, Wig, Van}).  The proof is standard, follows the scheme of Appendix 2 and we omit it.

The semiflow $S^t$ restricted to ${\mathcal M}$ is defined by the equations
\begin{equation}
\frac{dY}{d\tau}=   \xi F(Y,  \xi),
\label{inert}
\end{equation}
where
$$
    F(Y,  \xi)={ K}^{(1)} (Y) + \xi { K}^{(2)} (Y, \tilde { K}^{(1)}(Y)  +  W(Y, \xi)) +
$$
$$
    +
   \xi^2 { K}^{(3)}(\tilde { K}^{(1)}(Y)  +  W(Y,\xi)) + { R} Y + { T}\tilde { K}^{(1)}(Y)  + W(Y, \xi))
 +  f.
$$
The estimates for $W$ show that $F$ can be presented as
\begin{equation}
F(Y,  \xi)={ K}^{(1)} (Y)  + { R} Y + { T}\tilde { K}^{(1)}(Y) + f + \phi_{\xi}(Y)
\label{inert2}
\end{equation}
where a small correction $\phi_{\xi}$ satisfies
\begin{equation}
|\phi_{\xi}|_{C^1({\mathcal B}^p(R_0))} < c_0\xi^{1/2}.
\label{inert2c}
\end{equation}
In (\ref{inert2})  $ R$ and $f$ are free parameters. The quadratic form  ${ D}(Y)={ K}^{(1)} + { T}\tilde { K}^{(1)}$
 can be also considered as  a free parameter according to $p$- Decomposition Condition.
Therefore, we have proved the following assertion.

\begin{lemma} \label{quadr2} { Let
\begin{equation} \label{QField}
F(Y)= { D}(Y) + { R} Y + f
\end{equation}
be a quadratic vector field on ${\mathcal B}^p(R_0)$, where
$$
{D}_i(Y)=\sum_{j=1}^p \sum_{l=1}^p  D_{ijl} Y_j Y_l, \quad ({ R} Y)_i=\sum_{j=1}^p R_{ij} Y_j.
$$
Consider system  (\ref{gks2}), (\ref{gks3}). Let $p$- Decomposition Condition hold.
Then for any $\epsilon >0$ the field $F$ can be $\epsilon$ - realized by the semiflow $S^t({\mathcal P})$ defined  by
system (\ref{gks2}), (\ref{gks3}), where parameters $\mathcal P$ are the matrices ${ P}$, $ R$, $ \tilde  P$,
$ \tilde  R$ and the vectors $f, \tilde f$.
}
\end{lemma}

By  Lemma \ref{quadr2}  and results \cite{Stud} we prove the following assertion on realization of  all vector fields
by quadratic systems.

\begin{proposition} \label{QuadrP} { Consider the semiflows defined by  systems  (\ref{gks}),   where the triple $\{N, M, g \}$ serves as a parameter ${\mathcal P}$,  
for each $N$ the coefficients $K_{ijl}$ with $i, j, l \in \{1,..., N\}$ satisfy $p$-decomposition condition for a $p$ such that $N/2 < p^2 + p \le N$.  

Then these semiflows  enjoy the following property. 
For each integer $n$, each $\epsilon > 0$  and each vector field
$Q$ satisfying (\ref{cond1}) and (\ref{inward}), there exists a value of the parameter
$\mathcal P$ such that
the corresponding system   (\ref{gks}) defines a  semiflow $S^t_{\mathcal P}$, which $\epsilon$-realizes
the vector field $Q$ on $n$-dimensional positively invariant  manifold ${\mathcal M}_{n, Q}$.
 } 
\end{proposition}

Below we apply this result to the semiflows
defined by problem (\ref{OBEstream1})-(\ref{Maran3}).

\section{Control of linear terms  in system (\ref{maineq1})}
\label{sec: 10}

In this section we first show that the matrices ${\bf M}^{{\pm} {\pm}}$ involved in system  (\ref{maineq1})
are completely controllable by the function $u_1(x,y)$ and, thus, the linear terms  in the right hand side of
 this system satisfy the LOD condition from Sect. 4.

\subsection{\bf Control of matrix   ${\bf M}$ by $u_1$}
\label{sec: 10.1}

To calculate the entries of $M_{ij}$ 
  we take into account that this matrix
 can be decomposed to 4 blocks $M_{ij}^{\pm \pm}$,  where
\begin{equation}
   M_{ij}^{\pm \pm}(u_1(\cdot,\cdot))=
\langle \{\psi_j^{\pm}, \tilde \theta_i^{\pm}  \},  u_1\rangle.
\label{M}
\end{equation}
 A  straight forward calculation by (\ref{eigpms}), (\ref{eigpm})  and  (\ref{eigplus}) gives
\begin{equation}
   M_{ij}^{++}(u_1)=a_{ij}^{++}\int_0^{2\pi} \int_0^h
[\tilde \zeta_{i j}(y)\cos((i +j)x)
+  \zeta_{i j}\cos((i - j)x)]
 u_1(x,y)dx dy,
\label{M++}
\end{equation}
\begin{equation}
   M_{ij}^{--}(u_1)=a_{ij}^{--}
   \int_0^{2\pi} \int_0^h
[\tilde \zeta_{i j}(y)\cos((i + j)x)
  - \zeta_{i j}(y)
\cos((i - j)x)]
 u_1(x,y)dx dy,
\label{M--}
\end{equation}
and
\begin{equation}
   M_{ij}^{+-}(u_1)=a_{ij}^{+-}\int_0^{2\pi} \int_0^h
[\tilde \zeta_{i j}(y)\sin((i + j)x)
  +
\zeta_{i j}(y)\sin((i - j)x)]
 u_1(x,y)dx dy,
\label{M+-}
\end{equation}
\begin{equation}
   M_{ij}^{-+}(u_1)=a_{ij}^{-+}\int_0^{2\pi} \int_0^h
[-\tilde \zeta_{i j}(y)\sin((i + j)x)
 +
\zeta_{i j}(y)\sin((i - j)x)]
 u_1(x,y)dx dy,
\label{M-+}
\end{equation}
where
$$
\zeta_{ i j}=j\Psi_{j}(y)\frac{d\bar \theta_{i}(y)}{dy} + i \frac{d\Psi_{j}(y)}{dy} \bar \theta_{i}(y),
$$
$$
\tilde \zeta_{i j}=j\Psi_{j}(y)\frac{d\bar \theta_{i}(y)}{dy} - i \frac{d\Psi_{j}(y)}{dy} \bar \theta_{i}(y),
$$
and $a_{ij}^{\pm \pm}$ are  coefficients such that $|a_{ij}^{\pm \pm}|=1/2$. 
Using relations (\ref{eig3m}) and (\ref{asp}) for $\bar \theta_i$  and $\Psi_j$, one obtains
\begin{equation}
   \tilde \zeta_{i j}= \bar a_i \bar b_j \tilde \eta_{ij}, \quad \zeta_{i j}= \bar a_i \bar b_j \eta_{ij},
\label{M2}
\end{equation}
where $\bar a_i, \bar b_j$ are some non-zero coefficients, and
\begin{equation} 
    \eta_{i j}=ijy(\sinh(i+j)(h-y)) + \frac{hij \cosh(hi-(i+j)y)}{\sinh(jh)} - i\cosh(j(h-y))\cosh(i(h-y)),
\label{eta1}
\end{equation}
\begin{equation}
\tilde \eta_{i j}=ijy(\sinh(i-j)(h-y)) - \frac{hij \cosh(hi-(i-j)y)}{\sinh(jh)}  + i\cosh(j(h-y))\cosh(i(h-y)).
\label{eta2}  
\end{equation}

\begin{lemma} \label{8.1}
{
Given a quadruple of  $N \times N$ matrices $T^{\pm \pm}$ with  entries $T_{ij}^{\pm \pm}$ there
exists a $2\pi$ -periodic in $x$, smooth function $u_1(x,y)$
 such that the support of $u_1$ lies in the strip $\delta_1 < y  < h$ and
\begin{equation}
  M_{jl}^{+ +} [u_1] = T_{jl}^{++}, \quad  M_{jl}^{--} [u_1] = T_{jl}^{--},
\label{contM1}
\end{equation}
\begin{equation}
  M_{jl}^{+ -} [u_1] = T_{jl}^{+ -},  \quad M_{jl}^{-+} [u_1] = T_{jl}^{-+},
\label{contM2}
\end{equation}
 where $j, l =1,2,..., N$.
}
\end{lemma}

{\bf Proof}.
Let us show that systems of equations   (\ref{contM1}) and (\ref{contM2})
are resolvable.
 These two systems are independent, and we consider
the first one,  the  second one can be treated in a similar way.
 We represent $u_1(x,y)$ by  a Fourier series:
$$
u_1(x,y)=u_0^+(y) + \sum_{k=1}^{+\infty} \hat u_{k}^{+}(y) \cos( k x) + \hat u_{k}^{-}(y) \sin( k x).
$$
 By elementary transformations,
we reduce (\ref{contM1}) to the following form
\begin{equation}
 \int_0^h
\tilde \eta_{i j}(y) \hat u_{i+j}^+(y) dy=A_{ij}, \quad 1 \le i,j \le N,
\label{M11}
\end{equation}
\begin{equation}
 \int_0^h
 \eta_{i j}(y) \hat u_{j-i}^+(y) dy=B_{ij}, \quad N \ge j \ge i \ge 1,
\label{M21}
\end{equation}
\begin{equation}
 \int_0^h
 \eta_{i j}(y) \hat u_{i-j}^+(y) dy=C_{ij},  \quad N \ge i > j \ge 1,
\label{M31}
\end{equation}
where ${\bf A}, \ {\bf B}$ and ${\bf C}$ are some matrices.
In (\ref{M11}) we introduce an index $m$ by $m=i+j$, in (\ref{M21}) by $m=j-i$ and in (\ref{M31}) by $m=i-j$. These equations
become 
\begin{equation}
 \int_0^h
\tilde \eta_{i, m-i}(y) \hat u_{m}^+(y) dy=A_{i, m-i}, \quad 1 \le i < m,
\label{M11a}
\end{equation}
where $m \in \{2,...,2N\}$,
\begin{equation}
 \int_0^h
 \eta_{i, i+m}(y) \hat u_{m}^+(y) dy=B_{i, m+i}, \quad   1 \le i \le  N-m,
\label{M21a}
\end{equation}
where $m \in \{0, 1,..., N-1\}$,
\begin{equation}
 \int_0^h
 \eta_{i, i-m}(y) \hat u_{m}^+(y) dy=C_{i, i-m},  \quad   m < i \le N,
\label{M31a}
\end{equation}
where $m=0,1,..., N-1$.

Let us show that system (\ref{M11a})-(\ref{M31a}) has a solution using
the Fredholm alternative. The left hand sides of  (\ref{M11a})-(\ref{M31a}) define a linear map $U$ from on the set of $C^1$-smooth
functions $\hat u_m(y)$ defined on the interval $\delta_1 \le y \le h$
to a finite dimensional linear Euclidian space $\mathcal E$.
Let us consider the image $R(U)$ of $U$. If the closure $Clos(R(U))$
does not coincide with  the whole $\mathcal E$, then there is a nonzero vector from ${\mathcal E}$ orthogonal to $R(U)$. Therefore, in this case
there are coefficients $X_{m,i}, Y_{m,l}$ and  $Z_{m,k}$ such that
\begin{equation}
 \sum_{i=1}^{m-1} |X_{m, i}| + \sum_{i=1}^{N-m} |Y_{m, i}| + \sum_{i=m+1}^{N} |Z_{m, i}|=1,
\label{Zmi}
\end{equation}
\begin{equation}
S_m(y)  + \tilde S_m(y) \equiv 0
\label{Fred}
\end{equation}
for  all $y \in (\delta_1, h)$ and  all $m=\{1,..., N\}$,
where
\begin{equation}
  \tilde S_m=\sum_{ 1 \le i < m}  X_{m, i}\tilde \eta_{i, m-i}(y),
 \label{tS5}
\end{equation}
\begin{equation}
  S_m= \sum_{ 1 \le i < N-m}  Z_{m, i} \eta_{i, m+i}(y)
  + \sum_{ m < i \le N}  Y_{m, i} \eta_{i, i-m}(y).
\label{S5}
\end{equation}
Since  functions $S_m, \tilde S_m$  are analytic in $y$,
 equation (\ref{Fred}) is fulfilled  for all $y \in (-\infty, \infty)$.
Therefore, eq. (\ref{Fred}) means  that  nontrivial linear combinations of the functions
$\eta_{i, m+i}(y)$, $\tilde \eta_{i, m-i}(y)$ and $\eta_{i, i-m}(y)$ with coefficients $Z_{m,i}, X_{m,i}$ and $Y_{m,i}$
are zero  for all $y$ and $m$.

Let us prove that these nontrivial linear combinations do not exist.
We see  that
$$
\eta_{i, m+i}=\frac{h(i+m)i}{\sinh((m+i)h))}\cosh(hi -(2i+m)y) + i(i+m)y \sinh(2i+m)(h-y)) -$$
$$
-i\cosh((m+i)(h-y)) \cosh i(h-y)
$$
where $m+ 2i \in I_{m, N}=\{m+2, m+4, ... , 2N-m \}$,
$$
\eta_{i, i-m}=\frac{h(i-m)i}{\sinh((i-m)h))}\cosh(hi -(2i-m)y) + i(i-m)y \sinh(2i-m)(h-y)) -
$$
$$
-i\cosh((i-m)(h-y)) \cosh (i(h-y))
$$
where $2i-m \in  \{m+2, m+4, ... , 2N -m \} $,
and
$$
\tilde \eta_{i, m-i}=-\frac{h(m-i)i}{\sinh((m-i)h))}\cosh(hi +(m-2i)y) + i(m-i)y \sinh(2i-m)(h-y)) +
$$
$$
+i\cosh((m-i)(h-y)) \cosh (i(h-y)),
$$
where $m-2i \in J_m=\{m-2, m-4, ..., -m+2 \}$.
Let us observe that
functions $\cosh(a+ k_1y), \cosh(b+ k_2 y), \sinh(a' + k_1 y), \sinh(b'+k_2 y)$ and $ y\sinh(c+ k_1 y)$ are linearly independent
for all $a, b, a', b', c$, if $|k_1| \ne |k_2|$. The sets $I_{m, N}$ and $J_m$ are disjoint.
Thus the functions $S_m$ and $\tilde S_m$ are mutually linearly independent for any choice of
coefficients $X_{ij}, Y_{ij}$ and $Z_{ij}$ such that  $\sum_{ 1 \le i < m}  X_{m, i} >0$,
and $\sum_{ 1 \le i < N-m}  |Z_{m, i}| + \sum_{ m < i \le N}  |Y_{m, i}| >0$. We conclude thus that
either
\begin{equation}
\tilde S_m(y)\equiv 0, \quad y \in (-\infty, \infty),
\label {Fred1}
\end{equation}
or
\begin{equation}
S_m(y)\equiv 0 \quad y \in (-\infty, \infty).
\label {Fred2}
\end{equation}

Let us consider the case when (\ref{Fred1}) holds. Consider the function $\tilde S_m$ and  terms proportional  to $y$ in this function. Suppose that
there is a coefficient $X_{m, i} \ne 0$. We notice then, taking into account  only the terms $i(m-i)y \sinh(2i-m)(h-y))$
that (\ref{Fred1}) entails
$$
X_{m, i}= X_{m, m-i}.
$$
This relation implies, by the aforementioned   linear independency of hyperbolic functions, that the identities
 \begin{equation}
  X_{m,i} w_{i, m}(y)=0,  \quad for \ all \  y \in {\bf R},
\label{XF}
\end{equation}
hold for each $i \in \{ \frac{m}{2}-1, ..., m-1 \}$ and $m$.
Here
$$
w_{i,m}(y)=-\frac{h(m-i)i}{\sinh((m-i)h)}\cosh(hi +(m-2i)y)-
$$
$$
- \frac{h(m-i)i}{\sinh(ih)}\cosh(hi+(i-2m)y)
+ m\cosh((m-i)(h-y))\cosh(i(h-y))
$$
 Let us show that the function $w_{i,m}(y) \ne 0$ for some $y$.
This function can be represented as a sum of exponents
$\exp((m-2i)y)$, $\exp(-(m-2i)y)$, $\exp(my)$ and $\exp(-my)$ with some coefficients. The coefficient
at $\exp(my)$ is not zero, therefore, $w_{i,m}(y)$ is not zero for some $y$. Then eq.(\ref{XF})
implies that all coefficients $X_{m, i}=0$.

Let us consider now the case when relation (\ref{Fred2}) holds. We consider in $S_m(y)$  terms proportional to $y$. This gives
$
Y_{m, i+m}=-Z_{m, i}.
$
This relation entails
\begin{equation}
Z_{m, i}v_{i, m}(y)=0,  \quad for  \ all  \ y \in {\bf R},
\label{ZF}
\end{equation}
that hold for all  $i =1,2,...,m$ and $m$.
Here
$$
v_{i, m}(y)=h(m+i)i (\frac{1}{\sinh i h} - \frac{1}{\sinh (i+m)h })\cosh(hi - (2i +m)y) -
$$
$$
- m \cosh(2i+m) (h-y)).
$$
It is clear
that $v_{i,m}(y)$ is not zero for some $y$.
Then eq. (\ref{ZF}) entails $Z_{m, i}=0$ for all $i, m$.
The proof of the lemma is complete.

Let us formulate a lemma about control $f$ by $\eta_1$.

\begin{lemma} \label{fcontrol}
{ Given vectors $f^{+}=(f^{+}_1, ..., f_{N}^{+})$ and  $f^{-}=(f^{-}_1, ..., f_{N}^{-})$, there exists a
 smooth $2\pi$-periodic in $x$  function $\eta_1(x, y)$ with the support in the domain $\{(x,y): x \in (-\infty, +\infty),
 \delta_0 < y < h \}$, where $\delta_0 \in (0, h/2)$, such that
$$
\langle  \tilde \theta_i^{\pm}, \eta_1 \rangle =  f_i^{\pm}, \quad i=1,..., N.
$$
}
\end{lemma}
We omit an elementary proof.

\subsection{\bf Verification of $p$-Decomposition condition for system (\ref{G+}), (\ref{G-})}
\label{sec: 10.2}

To complete the proof of main result, it is necessary  to verify $p$-Decomposition condition from Section 4.
We choose the set $I_p$ and $i_l$ with $l=1,..., p$ in $p$-Decomposition condition by
$
i_l=k_l, \quad   I_p=\{k_1, ..., k_p\},
$
where  $k_i \in \{1,..., N\}$.  Let us set
$
Y_l= X^{+}_{k_l}.
$
Respectively, all the rest variables $X_i^{-}$ and $X^{+}_j$ with $j \ne k_l$ will be $Z_l$.
 Let us verify relation (\ref{barK}).
Note that the matrix $\tilde { K}^{(1)}$ involves $G_{jli}^{+++}$ and $G_{jli}^{+--}$.  
Therefore  it is sufficient to  verify that the linear system 
\begin{equation} \label{eqG}
\sum_{i \in J_p} G_{ijl}^{+++} u_i  =b_{jl}, \quad j, l \in I_p
\end{equation}
 has a solution for any given $b_{jl}$.
To check it,
let us calculate
the coefficients $G_{ijl}^{+++}$ defined by  (\ref{G+}).
 Integrating by parts one has
\begin{equation}
   G_{ijl}^{+++}=-\langle \{\psi_j^+, \tilde \theta_i^+ \}, \theta_l^+ \rangle + O(\nu^{-1}),
\label{GM11}
\end{equation}
thus  by definition (\ref{MatrM}) of $M_{ij}^{\pm \pm}$ one has
\begin{equation}
   G_{ijl}^{+++}=-M_{ij}^{++}( \theta_l^+)  + O(\nu^{-1}).
\label{GM12}
\end{equation}
Using this relation and (\ref{M++})
one obtains
$$
   2G_{ijl}^{+++}=   \int_0^{2\pi} \int_0^h
[\tilde \zeta_{ ij}(y)\cos((j + i)x)
  + \zeta_{ ij}(y)
\cos((j - i)x)] \Theta_l(y)\cos (lx)
 dx dy
+O(\nu^{-s}),
$$
where $s >0$.
We can transform this relation to the form
\begin{equation}
   G_{ijl}^{+++}=\frac{1}{4}(
       \delta_{l, i+j}  \tilde I_{ijl}
   +  \delta_{l, i- j}   I_{ijl} +  \delta_{l, j-i}  I_{ijl}) + O(\nu^{-s_1}), \quad s_1 >0,
\label{GM++2}
\end{equation}
where
 $\delta_{i,j}$ denotes the Kronecker symbol
and
$$
 \tilde I_{ijl} =\int_0^h  \tilde \zeta_{ij}(y)  \Theta_{l}(y) dy,   \quad  I_{ijl}= \int_0^h  \zeta_{ij}(y)  \Theta_{l}(y)   dy, \quad 1 \le i, j \le N
$$
Let us estimate   $\tilde I_{ij,i+j}$ and
$I_{ijl}$ for $l=i-j$ and $l=j-i$. 
We compute these integrals   taking into account $\nu >>1, h=\log \nu$ in relations (\ref{M2})
, (\ref{eta1}) and (\ref{eta2}) for
for $\xi_{ij},\tilde \xi_{ij}$.  As $\nu \to +\infty$ we have  the asymptotics
$$
\xi_{ij}= a_{i} b_j (2ijy -i)( 1+ O(\nu^{-s})),
$$
$$
\tilde \xi_{ij}=  a_{i}  b_j i (1 + O(\nu^{-s})), 
$$
where constants $a_{i}, b_j$ do not depend on $\nu$  and $s>0$.
For $\Theta_j$ we use  relation
 (\ref{asttheta2}).
We obtain then
\begin{equation}
 \tilde I_{ij, i+j}=a_i b_j  (\frac{i}{2(j+i)} + O(\nu^{-s})), \quad \nu  \to +\infty,
 \label{Iij}
\end{equation}
\begin{equation}
  I_{ij, i-j}= a_i b_j (\frac{j-i}{2i} + O(\nu^{-s})), \quad \nu  \to +\infty, \quad i \ge j
 \label{tIij}
\end{equation}
and
\begin{equation}
  I_{ij, j-i}= O(\nu^{-s}),  \quad \nu  \to +\infty, \quad  j \ge i.
 \label{Jij}
\end{equation}

Let us apply now Lemma  \ref{quadr2}. For each integer $p$ and a large $N=N(p)$
we solve  system (\ref{eqG}).  By (\ref{GM++2}) this system can be reduced to
the following form:
\begin{equation} \label{XX1}
   u_{k_i + k_j} \tilde I_{k_i, k_j, k_i+ k_j} + u_{k_i-k_j} I_{k_i, k_j, k_i -k_j} = b_{ij}, \quad
       i, j=1,..., p,  k_i > k_j
\end{equation}
\begin{equation} \label{XX2} 
   u_{k_i + k_j} \tilde I_{k_i, k_j, k_i+ k_j} + u_{k_j-k_i} I_{k_i, k_j, k_j -k_i} = b_{ij}, \quad
      i, j=1,..., p,  k_j > k_i
\end{equation}
\begin{equation} \label{XX3} 
   u_{2k_i } \tilde I_{k_i, k_i, 2k_i}  = b_{ii}, \quad
      i=j=1,..., p.
\end{equation}
We use a lemma.

\begin{lemma} \label{8.3}
{ There exists a subset ${\mathcal K}_p$ of $\{1,..., N\}$ satisfying the following conditions:

({\bf Ki})
$$
k_1 < k_2 < ... < k_p,  \quad      k_i=1 \quad  mod \ 3,  
$$

({\bf Kii}) all the sums $k_i +k_j$ are mutually distinct, i.e.,
$
k_i + k_j=k_{i'} + k_{j'}
$
implies $i=i', j=j'$  (in other words,  the map $(i,j) \to k_i +k_j$ from  $\{1,..., p\} \times \{1,..., p\}$ to $\{1,..., 2N\}$ is injective).

}
\end{lemma}

{\bf Proof}.
The set ${\mathcal K}_p=\{k_1, k_2, ..., k_p\}$ of  indices satisfying these conditions can be found by an iterative procedure.
Let us set $k_1=1, k_2=7$. We take a sufficiently large $k_3=3m_3+1$
such that  $3k_2 < k_3$  (for example, $k_3=19$). At $j+1$-th step we take an odd $k_{j}$  such that
$k_{j+1} > 3k_{j}$.
Notice that at $j$-th step all possible sums form a subset  of the set $\{2, ..., 2k_{j}\}$ and all possible differences give
a subset of the set $\{1, ..., k_{j}\}$.  Therefore, the new sums $k_{j+1}+ k_l$ do not coincide with the sums obtained earlier.
    This shows that
the inductive process can be continued, and completes the proof.

By this lemma we can prove resolvability of system (\ref{XX1})-(\ref{XX3}). 
Using {\bf Ki} we observe that $k_i + k_j \ne k_l - k_s$ for all $i, j,l,s \in \{1,..., p\}$. Due to this fact,  we can put  $u_{k_i - k_j}=0$ for all $i,j$ that simplifies system 
 (\ref{XX1})-(\ref{XX3}) and gives
\begin{equation} \label{KK1c}
   u_{k_i + k_j} \tilde I_{k_i, k_j, k_i+ k_j} = b_{ij}, \quad
       i, j=1,..., p.
\end{equation} 
By (\ref{Iij}),(\ref{tIij}) the coefficients $\tilde I_{k_i, k_j, k_i+ k_j}$
 are not zero for $i \ne j$.
Then due to {\bf Kii} 
equation (\ref{KK1c}) has a solution and thus $p$-Decomposition condition is fulfilled.

\section{Proof of main results}
\label{sec: 11}

Let us prove Theorem 3.2.
To make our arguments more transparent, we describe here a sequence of steps for $\epsilon$ -realization of
a vector field $Q$ (see (\ref{ordeq})) 
 on the unit ball
${\mathcal B}^n$. We suppose  $Q$ satisfy  (\ref{cond1}) and (\ref{inward}).
Our goal is to find parameters $\mathcal P$
of IBVP  (\ref{OBEstream1}) -(\ref{Maran3}) 
such that the coressponding family of  semiflows, generated by this IBVP,   
$\epsilon$ -realizes $Q$. We proceed it in two steps.
\vspace{0.2cm}

{\em Step 1}. 
Using proposition (\ref{QuadrP}), for  $\epsilon_0 > 0$ we $\epsilon_0$- realize the field $Q$  by  the family of the semiflows defined by quadratic vector fields    
(\ref{gks})  on ${\bf R}^N$ satisfying $p$-decomposition  condition  and with parameters $N, M, g$.   

Lemma \ref{8.3}  shows that systems (\ref{maineq1}) satisfy $p$-decomposition condition, thus, the corresponding semiflows  $\epsilon$-realize $Q$. 
The corresponding positively invariant manifold ${\mathcal M}_{n, Q}$ is diffeomorphic to the unit ball ${\mathcal B}^n$.
The right hand sides of system (\ref{maineq1})  is defined by a   quadratic field $V$,  which can be represented as $V(X)=K(X) + MX +g$,
(see (\ref{gks}) with some  $K=K(Q, \epsilon_0)$, $M=M(Q, \epsilon_0)$  and $g=g(Q, \epsilon_0)$.

{\em Step 2}. Consider the field $V$ defined via $N$, $K$, $M$ and $g$ and obtained at the previous step. For each $V$ and $\epsilon_1>0$  we find a IBVP (\ref{OBEstream1}) -(\ref{Maran3}) such that the corresponding semiflow defined by this problem $\epsilon_1$-realizes  the vector field  $V$ on
the locally invariant in the domain ${\mathcal D}_{\gamma, R_0, C_1, C_2, \alpha} \subset \mathcal H$  and locally attracting manifold
${\mathcal M}_{2N+1, \gamma}$ of dimension $2N+1$.
To  proceed it, we use the following  parameters: functions $U(y), u_1(x, y)$ satisfying (\ref{suppU}), (\ref{suppu1}),$\eta(x,y)$  and the numbers $\nu, \gamma$.  Here we  use the method of Section \ref{sec:4}. We first choose a sufficiently large $\nu > \nu_0$ (in order to use
Proposition \ref{6.2} about the spectrum of the linear operator $L$). Then for sufficiently small positive  $\gamma$ such that  $\gamma < \gamma_0(R_0, 
 \nu, \epsilon_1, \alpha, N, K, M, g)$  we can reduce the semiflow defined by    the IBVP  (\ref{OBEstream1}) -(\ref{Maran3}) to   system  (\ref{maineq1}).  For sufficient small $\gamma$
we can remove corrections $\phi_i$ in the right hand sides of (\ref{maineq1}). 
Note that by variations of $u_1$ and $\eta$  we can obtain all quadratic fields $V=KX + MX +g$ with  any prescribed   
$M$ and $g$ that follows from results of Sect. \ref{sec: 10.1}  (see Lemmas \ref{8.1} and \ref{fcontrol}).

If $R_0$ is large enough, we obtain  the embeddings
$$
{\mathcal M}_{n, Q}  \subset {\mathcal M}_{2N+1, \gamma} \subset {\mathcal H},
$$
where $\mathcal H$ is the phase space defined by  (\ref{phasespace}).
For sufficiently small $\epsilon_0, \epsilon_1$ 
 steps 1 and 2  give $\epsilon$-realization of the field $Q$ by semiflows defined by     IBVP   (\ref{OBEstream1}) -(\ref{Maran3}) with parameters
 $U(\cdot), u_1(\cdot, \cdot), \eta(\cdot,\cdot), \nu, \gamma$.    Theorem \ref{maint} is proved.


{\bf Proof} of {\bf Corollary}.   Consider a  flow on finite dimensional  smooth compact manifold defined by a $C^1$-smooth vector field and having a  hyperbolic compact invariant set $\Gamma$.  
For some integer $n>0$ we can find a smooth vector field $Q$ on a unit ball ${\mathcal B}^n$, which generates a semiflow having a topologically equivalent hyperbolic compact invariant set $\Gamma'$
(and the corresponding restricted dynamics are orbitally topologically equivalent) .  
Due to the Theorem  on Persistence of Hyperbolic sets (see \cite{Ru, Katok})  
we find a sufficiently small $\epsilon(\Gamma', Q) >0$ such that for all 
$C^1$ perturbations  $\tilde Q$ of $Q$ satisfying $|\tilde Q|_{C^1({\mathcal B}^n)} <\epsilon$ the perturbed vector fields $Q +\tilde Q$ generate
hyperbolic compact invariant sets $ \tilde \Gamma$   topologically equivalent  to $\Gamma$ (and the corresponding restricted dynamics are orbitally topologically equivalent). Then we $\epsilon$- realize this field by Theorem \ref{maint}.


\section{Conclusion}
\label{sec: 12}

In the case of a free liquid surface in contact with air, buoyancy and surface tension effects interplay in the convection. In this paper we have considered the model, where buoyancy is zero, then we are dealing with  the Marangoni effect and Marangoni -B\'enard convection studied in a number of physical and mathematical
works. 
It is shown that
the corresponding Navier -Stokes equations can generate  semiflows with complicated hyperbolic dynamics.
They can realize, with arbitrary accuracy, any finite dimensional vector fields.  This realization can  be done
on  stable invariant manifolds (in general, these manifolds are not globally attracting).
The main instrument in this realization is a choice of  an external heat source and the Prandtl number. The author thinks that the   methods  of this paper work for  
 the Boussinesq model with no-slip or free surface boundary conditions, but  it  is not clear how to apply the method of this paper for incompressible viscous fluids. 
Note that, in  Sections \ref{sec:4} and \ref{sec:5}, a general method to prove existence of chaotic behaviour for quadratic systems of ODE's and PDE's is stated.

\section{Acknowledgements}

I dedicate this paper to  memory of my friends Vladimir Shelkovich and Sergey Aristov.
I am  grateful to Prof. D. Grigoriev (Lille) for the help and discussions and A. Weber (Bonn) for interesting discussions.
Also the author is thankful to Prof. V. Kozlov (Linkoping) for critics (severe, but useful), and
 Prof. V. Volpert (Lyon) for his hospitality in Lyon University during 1998-2003 and fruitful discussions.
Also I am grateful the anonymous Referees for useful remarks.

 I am grateful to Prof. D. Ruelle, who has encouraged, in 1991,  the first steps of the  author in this domain.

This work  was financially supported by Government of Russian Federation, Grant 074-U01.


\section{Appendix 1 }

Let us consider the Fourier series for $\bar \omega$ and $\bar \psi$. Let us represent  $u$  and $Tr(u)$ by
\begin{equation}
u=\sum_{k\in {\bf Z}} \sum_{m=0}^{\infty} \hat u_{k, m} \cos (mh^{-1} y) \exp(i kx),
 \label{Ap1}
 \end{equation}
\begin{equation}
Tr(u)= \sum_{m=0}^{\infty} \hat u_{k} \exp(i kx), \quad \hat u_k=\sum_{m=0}^{\infty} \hat u_{k, m}.
 \label{Ap2}
 \end{equation}
 Then the vorticity $\bar \omega$, defined by (\ref{Stokes}), (\ref{omegab}), is
 \begin{equation}
\bar \omega= \sum_{k \in {\bf Z}} ik \hat u_{k} \exp(i kx) \frac{\sinh( k(h-y))}{\sinh kh}.
 \label{Ap3}
 \end{equation}
 We observe that
 \begin{equation}
||\bar \omega(u)||^2 \le C_1 \sum_{k \in {\bf Z}} |k| |\hat u_{k}|^2,
 \label{Ap4}
 \end{equation}
 for some $C_1> 0$.
 By (\ref{Ap3}) one obtains that the stream function $\bar \psi$ has the form
 \begin{equation}
\bar \psi=  \sum_{k \in {\bf Z}} i \hat u_{k} \exp(i kx) (-\frac{y\cosh( k(h-y))}{2\sinh kh} + \frac{h\sinh ky}{2(\sinh kh)^2}).
 \label{Ap5}
 \end{equation}
This relation implies
$$
||\nabla \bar \psi||^2 \le c_1 \sum_{k \in {\bf Z}} |k| |\hat u_k|^2  \le c_2 ||u||_{\alpha}, \quad \alpha > 1,
$$
(the second inequality follows from  embedding (\ref{Tr}) for traces).
Since $$
\sup_{y \in (0, h)} |ky| \cosh(k(h-y)) (\sinh kh)^{-1} \le C_3,
$$
one has
$$
|\nabla \bar \psi| \le C_1  \sum_{k \in {\bf Z}, k \ne 0}| \hat u_{k}| \le
(\sum_{k \in {\bf Z}, k \ne 0}|\hat u_{k}|^2 k^{2\gamma})^{1/2}  (\sum_{k \in {\bf Z}, k \ne 0} k^{-2\gamma})^{1/2}
$$
with $\gamma >1/2$.
The last estimate gives
$$
|\nabla \bar \psi| \le C_2 ||u||_{\alpha}, \quad  \alpha > 3/2.
$$

\section{Appendix 2}

{\bf Proof of Lemma \ref{rvf}}.  This assertion is a consequence of Theorem 6.1.7  \cite{He}.
In the variables $\tilde w$, $X$ system (\ref{eveq5}), (\ref{eveq6}) takes the form
\begin{equation}
X_t= g(X, \tilde w),   
\label{eveq5a}
\end{equation}
where 
$$
g=\gamma^{-1} {\bf P}_1 F(\gamma (X + w_0) + \tilde w) + \gamma \bar f_1,
$$
\begin{equation}
\tilde w_t=L \tilde w + f_0(X, \tilde w), \quad f_0= {\bf P}_2 F(\gamma (X + w_0) + \tilde w).
\label{eveq6a}
\end{equation}
Using a standard truncation trick we first modify  eq. (\ref{eveq5a}) as follows:  
\begin{equation}
X_t= G(x, \tilde w),    
\label{eveq5at}
\end{equation}
where
$$
G(x, \tilde w)=g(X, \tilde w)\chi_{R_1}(X)
$$
and $\chi_{R_1}$ is a smooth function such that $\chi_{R_1}(X) =1 $ for $|X| < R_1$ and  $\chi_{R_1}(X) =0 $ for $|X| >  2R_1$.
After this modification, $X$-trajectories of (\ref{eveq5at}) are defined for all $t \in (-\infty, +\infty)$ (as in Theorem 6.1.7 \cite{He}).
Then an invariant   manifold for the semiflow defined by system  (\ref {eveq5at}), \ref{eveq6a}) is a locally invariant one for the semiflow generated by
  (\ref {eveq5a}), \ref{eveq6a}). 

Let us consider the semigroup $\exp(Lt)$.
If $w(0) \in B_2$ we have the following estimates \cite{He}
$$
||\exp(Lt) w(0)|| \le C_0||w(0)|| \exp(-\beta t),
$$
$$
||\exp(Lt) w(0)||_{\alpha} \le C_1||w(0)|| t^{-\alpha}\exp(-\beta t),
$$
where $M, \beta >0$ do not depend on $\gamma$.
Moreover,
\begin{equation}
M_0=\sup_{(X, \tilde w) \in {\mathcal B}_{R, \gamma, C}}  ||f|| <  c_2\gamma^2,
\label{Ppr2}
\end{equation}
\begin{equation}
\lambda=\sup_{(X, \tilde w) \in {\mathcal B}_{R, \gamma, C}}  ||D_X f_0|| + ||D_{\tilde w} f_0|| <  c_3\gamma^2,
\label{Ppr3}
\end{equation}
\begin{equation}
M_2=\sup_{(X, \tilde w) \in {\mathcal B}_{R, \gamma, C}}  ||D_{\tilde w}  g|| <  c_4\gamma,
\label{Ppr4}
\end{equation}
We set   $\mu=\beta/4$. Then for small $\gamma$
\begin{equation}
M_3=\sup_{(X, \tilde w) \in {\mathcal B}_{R, \gamma, C}}  ||D_{X}  g|| <  c_4\gamma  < \mu.
\label{Ppr5}
\end{equation}
We set $\Delta=2\theta_1$, where
\begin{equation}
\theta_p =\lambda M_0 \int_0^{\infty} du u^{-\alpha} \exp(-(\beta - p\mu') u), \quad 1 \le p \le 1+r,
\label{Int1}
\end{equation}
and
$\mu'= \mu + \Delta M_2$. For sufficiently small $\gamma$ one has $\mu' <\beta/2$, therefore, the integral in the right hand side
of (\ref{Int1}) converges and, according to (\ref{Ppr4}), one obtains $\theta < c_5 \gamma^2$ (since $M$ is independent of $\gamma$).
We notice then that for sufficiently small $\gamma$ the following estimates
$$
(1+r) \mu' < \beta/2,
$$
$$
\theta_1 < \Delta(1 + \Delta)^{-1},   \quad \theta_1 < 1, \quad  \theta_1(1+\Delta)M_2{\mu'}^{-1} < 1,
$$
and
$$
\theta_p(1 +  \frac{(1+\Delta)M_2}{r \mu'}) < 1
$$
hold.


\end{document}